\documentclass{amsart}
\usepackage{amssymb,amscd,verbatim}
\usepackage{enumerate}

\newcommand{\pa}{\partial}

\newcommand{\CI}{\mathcal{C}^\infty}

\newcommand{\bS}{{}^b\kern-1pt S} 
\newcommand{\bT}{{}^b\kern-1pt T} 
\newcommand{\Hom}{\operatorname{Hom}}

\newcommand{\cB}{\mathcal{B}}

\newcommand{\cD}{\mathcal{D}}

\newcommand{\cH}{\mathcal{H}}

\newcommand{\cJ}{\mathcal{J}}
\newcommand{\cJi}{\cJ_{\infty}}
\newcommand{\cJs}{\cJ^{*}}
\newcommand{\cL}{\mathcal{L}}

\newcommand{\cP}{\mathcal{P}}

\newcommand{\fA}{\mathfrak{A}}

\newcommand{\hTa}[1]{\cH_{#1}(T)}
\newcommand{\hTk}{\hTa{k}}
\newcommand{\hTn}{\hTa{0}}

\newcommand{\hTu}{\hTa{\infty}}
\newcommand{\id}{{\rm id}}

\newcommand{\norm}[2]{\|#1\|_{#2}}

\newcommand{\otl}{\omega_{T}^{\ell}}
\newcommand{\otlr}{\omega_{T}^{\ell,r}}
\newcommand{\otr}{\omega_{T}^{r}}

\newcommand{\re}{{\rm Re}\,}

\newcommand{\sh}[1]{#1^{\sharp}}

\newcommand{\cal}{\mathcal}

\newcommand{\CC}{\mathbb C} 
\newcommand{\NN}{\mathbb N} 
\newcommand{\RR}{\mathbb R} 
 
\newcommand{\ZZ}{\mathbb Z} 
 
\newcommand{\be }{\begin{eqnarray*}} 
\newcommand{\ee }{\end{eqnarray*}}

\newcommand{\mB}{\mathfrak B}

\newcommand{\mH}{\mathcal H} 
\newcommand{\mL}{\mathcal L} 
\newcommand{\mP}{\mathcal P} 
\newcommand{\mR}{\mathcal R} 
\newcommand{\mW}{\mathcal W} 
\newcommand{\mO}{\mathcal O} 

\let\<\langle
\let\>\rangle

\let\rho\varrho

\newcommand{\ie}{{\em i.e., }}

\newtheorem{theorem}{Theorem}[section] 
\newtheorem{proposition}[theorem]{Proposition} 
\newtheorem{corollary}[theorem]{Corollary} 
 
\newtheorem{lemma}[theorem]{Lemma} 
 
\theoremstyle{definition}
\newtheorem{definition}[theorem]{Definition} 
\theoremstyle{remark} 
\newtheorem{remark}[theorem]{Remark} 
\newtheorem{example}{Example}

\long\def\komment#1{}

\def\be{{\beta}}

\def\phi{{\varphi}}

\def\pa{{\partial}}

\author[B. Ammann]{Bernd Ammann} \address{Universit\"at Hamburg, 
Fachbereich 11 -- Mathematik, Bundesstra\ss{}e 55, D - 20146 Hamburg, Germany,
{\rm http://www.berndammann.de/publications},}  
\email{ammann@math.uni-hamburg.de}
 
\author[R. Lauter]{Robert Lauter} \address{Universit\"at Mainz,
Fachbereich 17 -- Mathematik, D - 55099 Mainz, Germany}
\email{lauter@mathematik.uni-mainz.de}
        
\author[V. Nistor]{Victor Nistor} \address{Pennsylvania State
University, Math. Dept., University Park, PA 16802}
\email{nistor@math.psu.edu}

\author[A. Vasy]{Andr\'as Vasy} \address{M.I.T., Department of Mathematics,
Cambridge, MA
02139} \email{andras@math.mit.edu}
 
\thanks{V. N. was partially supported by NSF Grants DMS 99-1981 and
DMS 02-00808. A. V. was partially supported by NSF Grant DMS 02-01092
and by a Fellowship from the Sloan Foundation.
B.~A. was partially supported 
by The European Contract Human Potential Programme,
Research Training Networks HPRN-CT-2000-00101 and HPRN-CT-1999-00118.
}

\begin{document}


\title[Complex powers]{Complex powers and non-compact
manifolds}

\begin{abstract}  
We study the complex powers $A^{z}$ of an elliptic, strictly positive
pseudodifferential operator $A$ using an axiomatic method that
combines the approaches of Guillemin and Seeley. In particular, we
introduce a class of algebras, ``extended Weyl algebras,'' whose
definition was inspired by Guillemin's paper \cite{Guil85}. An
extended Weyl algebra can be thought of as an algebra of ``abstract
pseudodifferential operators.'' Many algebras of pseudodifferential
operators are extended Weyl algebras. Several results typical for
algebras of pseudodifferential operators 
(asymptotic completeness, construction of Sobolev spaces, boundedness 
between apropriate Sobolev
spaces, ... ) generalize to extended Weyl
algebras. Most important, our results may be used to
obtain precise estimates at infinity for $A^{z}$, when $A > 0 $ is
elliptic and defined on a non-compact manifold, provided that a
suitable ideal of regularizing operators is specified (a
submultiplicative $\Psi^*$--algebra). 
We shall use these results in a forthcoming paper to study
pseudodifferential operators and Sobolev spaces on manifolds with a
Lie structure at infinity (a certain class of non-compact manifolds).
\end{abstract}

\maketitle \tableofcontents

\section*{Introduction} 

Results about the complex powers $A^{z}$ of strictly positive,
elliptic pseudodifferential operators have proved useful to study the
asymptotic of eigenvalues of self-adjoint pseudodifferential operators
\cite{Guil85}, Sobolev spaces, and even certain non-linear
differential equations (see \cite{CSchrohe,Schrohe.Seiler,DoreVenni} and the
references therein).  The study of complex powers has been inaugurated
in Seeley's celebrated paper \cite{Seeley}, where he proved that if $A
> 0$ is an elliptic differential operator of order $m > 0$ on a
compact manifold, then $A^{z}$ is a pseudodifferential operator of
order $m \re(z)$. For previous results on complex powers, see
\cite{buci, Loya1, Loya2, Guil85, RS, SchroheC, ShubinBook, Seeley}.

In \cite{Guil85}, Guillemin has developed another approach to the
construction of complex powers. His approach is axiomatic, in the sense
that it works for operators in a ``Weyl algebra.'' Weyl algebras were
introduced in the aforementioned paper and are a generalization of
algebras of pseudodifferential operators on {\em compact}
manifolds. Thus, in particular, on compact manifolds, an operator of
negative order in a Weyl algebra is compact.

In this paper, we generalize Guillemin's approach to include
non-compact manifolds in which a certain algebra of regularizing
operators is specified. We thus drop the condition that an operator
of negative order be compact.  This forces us to replace the axioms of
a Weyl algebra, with the axioms of an ``extended Weyl algebra.'' In
particular, we have
to carefully check certain analytic facts that were
obvious in the case of a Weyl algebra.

An extended Weyl algebra should be thought of as an abstract algebra
of pseudodifferential operators. Many results on pseudodifferential
operators on a non-compact manifold generalize to extended Weyl
algebras. These algebras provide a natural framework to discuss
complex powers, Sobolev spaces, continuity, self-adjointness,
resolvents, holomorphic families of operators, ellipticity, and
parametrices for several classes of algebras of pseudodifferential
operators on non-compact manifolds. To be able to construct complex
powers, we also need to assume that our Weyl algebra satisfies two
additional properties, namely, that the kernel of the principal symbol map
consists of operators of lower order (a familiar property of
pseudodifferential operators) and that the given ideal of regularizing
operators is a submultiplicative Fr\'echet algebra that contains the
inverse of any $L^2$-invertible element (\ie that it is a
$\Psi^*$-algebra in the sense of Gramsch \cite{gra84, lmn2, Schrohe}).

An important role in our proofs is played by a theorem of Banach from
1948, \cite{ban48}, which states that inversion is continuous on a
Fr\'echet algebra whose set of invertible elements is a $G_{\delta}$
subset ($:=$ a countable intersection of open subsets). Any $\Psi^*$
algebra has an open set of invertible elements, so inversion is
continuous on a $\Psi^*$-algebra.

One of the main applications of our results will be to study the
complex powers of operators on manifolds with a Lie structure at
infinity \cite{alnv2}. The manifolds with a Lie structure at infinity
were introduced, informally, in \cite{meicm, MelroseScattering}, and
were studied more systematically in \cite{aln1}.

The paper is organized as follows. In the first section we introduce
the axioms of an extended Weyl algebra $\mW$, following \cite{Guil85},
and we include two examples. Anticipating, let us just say now that an
extended Weyl algebra is an algebra $\mW = \cup \mW^{\mu}$, $\mu \in
\CC$, of (possibly) unbounded operators with a common domain contained
in a given Hilbert space $\mH$ and such that $\mW^{\mu}\mW^{\nu} =
\mW^{\mu + \nu}$ and satisfying several axioms. An example is the
algebra of (sums of) classical pseudodifferential operators of complex
order on a compact manifold, with $\mW^{\mu}$ the space of operators
of order (at most) $\mu$. Certain algebras of pseudodifferential
operators on non-compact manifolds also provide us with examples of
extended Weyl algebras.

In Section \ref{Sec.FC} we discuss some of the basic properties of
extended Weyl algebras $\mW$, including boundedness in the given
Hilbert space. We prove in particular that any elliptic, symmetric
operator $T \in \mW^{m}$ (that is of order $m$), $m > 0$, is
automatically essentially self-adjoint. Intuitively, this means that
extended Weyl algebras model geometric operators on complete
manifolds. The setting of operators on cones (see \cite{Gil, Loya2,
MelroseAS}, for example), is not directly covered in our approach. The
domains of operators of positive order, which turn out to be some
generalizations of the usual Sobolev spaces, are studied in Section
\ref{Sec.AS}. 
In Section \ref{Sec.HFO}, we discuss holomorphic families
of operators, including asymptotic completeness. 
In Section \ref{Sec.SI} we proved that an elliptic
operator in $\mW^{\mu}$, $\re(\mu) \ge 0$, that is invertible as an
unbounded operator on the ambient Hilbert space, will have an inverse
in our extended Weyl algebra $\mW$, provided that the ideal of
regularizing operators is a spectrally invariant Fr\'echet algebra. We
also prove the continuity of the adjoint and the multiplication in
this section. In Sections
\ref{Sec.FC} -- \ref{Sec.HFO} we do not assume that multiplication
$S^m_{1,0} \times \mW^{-\infty} \ni (a, T) \to q(a)T$ is continuous
(this is Axiom \eqref{axiom.product}).  Holomorphic families of
variable order, called ``special holomorphic families,'' are studied
in Section \ref{Sec.SHF}. Special holomorphic families are then used
in Section \ref{Sec.CP} to construct complex powers of elliptic
operators, by combining the methods of \cite{Guil85, Seeley} and
\cite{ShubinBook}. In the process, we are led to study the behavior of
the resolvent $(T + \imath t)^{-1}$, where $T \in \mW^{m}$, $m > 0$,
is elliptic and symmetric, and we prove, in particular, that it is
bounded in $\mW^{-m}$, for $|t| \ge 1$. This allows us to recover the
results of \cite{Guil85, Seeley}.

It is worthwhile mentioning that the setting of \cite{Guil85,
Seeley} is such that their Weyl algebras are complete in the Guillemin
topology (\ie the topology on the set of regularizing operators is
given by the operator norms $H^{(s)} \to H^{(s')}$, for any $s,
s'$). This simplifies greatly the proofs. However, certain
applications, for example to manifolds with cylindrical ends
\cite{MitreaNistor} require us to go beyond that.  It would be
interesting and important to see to what extent the functional
calculus with symbolic functions of \cite{HV} extends to the setting
of extended Weyl algebras.

Throughout this paper, the notation ``$:=$'' means ``the left hand side
is defined to be equal to the right hand side.''


\section{Extended Weyl algebras\label{Sec.EWA}}

We introduce now a class of algebras generalizing the class of Weyl
algebras introduced in \cite{Guil85}. We follow the spirit of
\cite{Guil85} when defining these algebras. In fact, the main
difference is that we do not impose compactness conditions on
operators of negative order, which in turn forces us to make very
explicit analytic assumptions.  Then, in order to construct complex
powers, we need the ideal of operators of order $-\infty$ to satisfy
some strong analytic conditions that will be made precise below
(Conditions $(\sigma)$ and $(\psi)$).

\subsection{Notation and preliminaries}\
We need to recall first the definition of some classes of symbols.
Let $V \to M$ be a smooth orthogonal bundle over a
manifold $M$, possibly with corners, endowed with an orthogonal
connection $\nabla$. Let
\begin{equation}\label{def.<>}
	\<v\> = (1 + \|v\|^2)^{-1/2}\,, \quad v \in V.
\end{equation}
We shall denote by $S_{1, 0}^m(V)$ the space of type $(1, 0)$--symbols
of order $m$ on $V$ \cite{hor3, Taylor2}.  
Recall that the space $S_{1, 0}^m(V)$ is defined as follows. If
$V = M \times \RR^N$, then $S_{1, 0}^m(V)$ consists of the smooth
functions $a : V \to \CC$ such that 
\begin{equation}\label{eq.def.10}
	\| \pa_{x}^{\alpha} \pa_{\xi}^{\beta} (\phi_i(x) a(x, \xi)) \| 
	\le C(\alpha,\beta) \,\< \xi \>^{m - |\beta|},
\end{equation}
for any $x \in M$ and $\xi \in \RR^N$. The topology is given by the best
constants $C(\alpha, \beta)$ in Equation \eqref{eq.def.10}.

To define the space $S_{1,0}(V)$ and its topology in general, we first choose
a locally finite covering $(U_i)$ of $M$ by open subsets on which $V$ is trivial.
Then choose trivializations $\psi_i$ of $V$ on $U_i$ and a smooth partition of 
unity $\phi_i$ subordinate to the cover $U_i$. Then a smooth complex-valued 
function $a$ on the total space of the bundle $V$ is an element of $S_{1, 0}^m(V)$ 
if $\phi_i(x) a(\xi,x)$ is a type $(1,0)$--symbol of order $m$ for all $i$. This definition is
independent of the choice of the covering, trivialization, and partition of unity
used. If $M$ is compact, the topology is obtained by considering a finite partition 
of unity and the best constants in the definitions of symbols. (This topology is
thus a quotient topology.)

We shall also need the case of a manifold $M$ of bounded geometry and $V = TM$.
Then we specify the families $U_i$ and $\psi_i$ to be given by normal coordinates
around each point and requiring
\begin{equation}\label{eq.def.101}
	\| \pa_{x}^{\alpha} \pa_{\xi}^{\beta} (\phi_i(x) a(x, \xi)) \| 
	\le C(\alpha,\beta) \,\< \xi \>^{m - |\beta|},
\end{equation}
for any multi-indices $\alpha \in \NN^{l}$ and $\beta \in \NN^{k}$, 
and any $i$. Note that the definition of $S_{1, 0}^m(V)$ does not depend on 
the choice of $U_i, \psi_i, \phi_i$. 


Since most of our symbols are defined on $V$, we shall usually omit 
$V$ from the notation, thus write $S^{m}_{1, 0} := S^{m}_{1, 0}(V)$.

Also, recall that a symbol $a \in S^{m}_{1, 0}$ is called {\em
homogeneous of order (at most) $\mu = m + it$}, $t \in \RR$, if there
exists $R > 0$ such that
$$
	a(rv) = r^{\mu}a(v)\,,\quad \text{for all }\; r \in
	(0,\infty)\,, \text{ and } v \in V\,,\; \text{ such that }\,
	|t|, \|v\| \ge R.
$$
Moreover, a symbol $a \in S^{m}_{1, 0}$ is called a {\em classical
symbol of order $\mu = m + it \in \CC$,} $t \in \RR$, if there exist
symbols $a_{j} \in S^{m - j}_{1, 0} $, homogeneous of order $\mu -j$
such that
$$
	a(v) - \sum_{j = 0}^{k-1} a_j(v) \in S^{m-k}_{1, 0}.
$$  
(If this is the case, we write $a \sim \sum_{j = 0}^\infty a_{j}$.)
The subspace of classical symbols of order (at most) $\mu$ will be
denoted $S_{cl}^{\mu}(V) \subset S_{1, 0}^{m}(V)$, as usual.  We shall
often write $S_{cl}^{\mu} := S_{cl}^{\mu}(V)$.

The following lemma is probably well known (see
\cite[Lemma18.1.10]{hor3} for the case $m=0$).

\begin{lemma}\label{lemma.symb}\ 
Let $m \ge 0$, $k \in \RR$, and $S^{m}_{1, 0} := S^{m}_{1, 0}(V)$, as
above. Also, denote
$$
	{Ell}^{m}_{\mu} := S^{m}_{1, 0} \cap \{ |a(v)| \ge \mu^{-1}\<
	v \>^m \text{ if } |v| \ge \mu \}.
$$
Then
\begin{equation}
 	S^{k}_{1, 0}(\CC) \times {Ell}^{m}_{\mu} \to S^{mk}_{1, 0},
 	\qquad (f, a) \mapsto f \circ a
\end{equation}
is a well-defined and continuous for any $\mu > 0$. A similar
statement holds for $S^k_{1,0}(\CC)$ replaced with $S^{k}(\RR)$ or
with $S^{k}([\epsilon,\infty))$, if one restricts to real valued
symbols in $S^{m}_{1, 0}$ or, respectively, to symbols that are $\ge
\epsilon$.
\end{lemma}

\begin{proof}\ 
We can assume that $V$ is trivial. 
Thus it is enough to verify the
condition of Equation \eqref{eq.def.10}. We shall do this by induction
on $|\alpha| + |\beta|$.  Let $|\alpha| = |\beta| = 0$. By the
definition of $S_{1,0}^k(\CC)$ we have $|f(a(v))| \le C_1 \< a(v)
\>^{k}$. If $k\ge 0$, then $\< a(v) \>^{k} \le C_2 \< v \>^{km}$ gives
the desired inequality
\begin{equation}\label{eq.fav}
	|f(a(v))| \le C_3 \< v \>^{km}.
\end{equation}
Here the $C_i$'s are constants. For $k<0$ we use the ``ellipticity
condition'' $|a(v)| \ge \mu^{-1}\< v \>^m$ to conclude \eqref{eq.fav}.

The inductive step follows from the identity
$$
	X (f \circ a) = (f' \circ a) X(a),
$$
for any vector field $X$ on $V$.

The proof of the last two statements of the lemma follow by
restriction.
\end{proof}


If $X$ and $Y$ are two locally convex spaces, we shall denote by
$\mL(X,Y)$ the space of continuous linear maps $X \to Y$.  We shall
write $\mL(X,X) = \mL(X)$. Also, we shall denote by $T^*$ the adjoint
of a (possibly unbounded) operator $T$ between Banach spaces.

\subsection{The axioms of an extended Weyl algebra}\
We now introduce the conditions defining an ``extended Weyl algebra.''

Recall that a limit-of-Frechet space or, for short, an 
{\em $LF$--space} $X = \cup_{k \ge 1} X_k$ is a 
inductive limit $X = \displaystyle{\lim_\to} X_k$, $k = 1, 2, \ldots$,
where each $X_k$ is a Fr\'echet space, $X_k \subset X_{k+1}$ is a
closed subspace whose topology coincides with is the subspace topology. 
We endow $X$ with the inductive limit topology.
This implies, in particular, that $X_k$ also carries the subspace 
topology of $X$, or equivalently by defintion, 
that $X$ is a \emph{strict} inductive limit. (See, for
example, \cite[II.6.1--6.6]{Schaeffer:Topological} for a thorough
discussion.) 
By an {\em $LF$-algebra} we shall mean an $LF$--space $A =
\cup_{n \ge 1} A_n$, which is endowed with an algebra structure such
that $A_{k}A_{l} \subset A_{k + l}$ and the multiplication $A_{k}
\times A_{l} \to A_{k + l}$ is continuous, or equivalently $A\times A
\to A$ is continuous.  \\

In this section we consider a tuple $({\cal W}^{-\infty},{\cal
H},V,{\cal P}, q)$ satisfying the following axioms:
\begin{enumerate}[(i)]
\item\ \label{axiom.H} $\mH$ is a Hilbert space and 
$\mW^{-\infty} = \cup_{k \ge 1} \mW^{-\infty}_{k}$ is an $LF$--algebra
continuously embedded in $\mL(\mH)$ such that 
${}^* : \mL(\mH)\to \mL(\mH)$, $T\mapsto T^*$,
maps $\mW^{-\infty}$ continuously onto itself.
%
%
\vspace*{2mm}
\item\ \label{axiom.R} There exists an injective operator $R = R^* \in
\mW^{-\infty}_1$.
\vspace*{2mm}
\item\ \label{axiom.q} $V \to M$ is a vector bundle over
a compact manifold $M$, possibly with corners, or $V = TM$ and $M$ is with 
bounded geometry. Furthermore, $q$ is a map $q : S_{1, 0}^\infty \to
\Hom(\mW^{-\infty}\mH, \mH)$, $q(1) = I$, satisfying:
\begin{enumerate}[(a)]
\item $q(a)$ is symmetric for $a$ real; and
\item The restriction $q\vert_{S^{-\infty} } : S^{-\infty} \to
\mW^{-\infty}_{1}$
is well defined and continuous.
\end{enumerate}
\vspace*{2mm}
\item \label{axiom.product'}\ $q(S^m_{1,0})\mW^{-\infty}_{k}
\subset \mW^{-\infty}_{k + 1}$.
\item \label{axiom.P} The map $\mP: S_{1,
0}^\infty \times S_{1, 0}^{\infty} \to S_{1, 0}^{\infty} $ is bilinear and satisfies:
\begin{enumerate}[(a)]
\item $q(a) q(b) - q(\mP(a,b)) : S_{1, 0}^\infty \times S_{1,
0}^{\infty} \to \mW^{-\infty}_{1}$ is well defined and continuous;
\item The induced maps 
	$\mP : S_{1, 0}^m \times S_{1, 0}^{n} \to S_{1, 0}^{m+n},$
	$\mP : S_{cl}^{\mu} \times S_{cl}^{\nu} \to S_{cl}^{\mu + \nu},$
	and  $\mP_0 : S_{1, 0}^m \times S_{1, 0}^{n} \to S_{1,
	0}^{m+n-1},$
$\mP_0(a,b) = \mP(a,b) - ab$, are well defined and continuous.
\end{enumerate}
\vspace*{2mm}
\item \label{axiom.cont} There exists $d>0$ and a continuous seminorm
$p$ on $S_{1, 0}^{-d-1}$ such that 
$q(a)$ extends to a bounded operator on $\mH$ with  
$\|q(a)\| \le p(a)$, for all $a \in
S_{1, 0}^{-d-1}$.
\vspace*{2mm}
\item \label{axiom.product} The induced map
$$
	S_{1, 0}^m \times \mW^{-\infty}_{k} \to \mW^{-\infty}_{k + 1}, \qquad 
       (a,T) \mapsto q(a)T,
$$ 
(which is well defined by axiom~\eqref{axiom.product'}) 
is continuous.
\end{enumerate}
\vspace*{3mm}

Note that Axiom~\eqref{axiom.R} implies that $\mW^{-\infty}\mH$ is
dense in $\mH$ because $R\mH$ is dense for any $R = R^*$ satisfying
Axiom~\eqref{axiom.R}. Hence, $q(a)$ is an operator densely defined in
$\mH$ and the following axioms
make sense.\\


\noindent {\bf Comments.}\ Let us observe that the
Axiom~\eqref{axiom.product'} implies that if $T \in \mW^{-\infty}_{k}$
and $a \in S^{\infty}_{cl} := S^{\infty}_{cl}(V)$, then $T q(a) =
(q(\overline{a})T^*)^* \in \mW^{-\infty}_{k + 1}$. In particular,
$\mW^{-\infty}$ is closed under both left and right multiplication by
operators of the form $q(a)$.\\

Two other conditions, Conditions $(\sigma)$ and $(\psi)$ will be
considered below and will be used to construct complex powers, but
they are not part of the definition of an extended Weyl algebra. Before
stating these conditions, let us notice that the axioms above give
rise to an algebra.

\begin{proposition}\label{prop.def.mW}\ 
Let $\mH$, $\mW^{-\infty} = \cup_{k \ge 1} \mW^{-\infty}_{k}$, and $q$
satisfy the Axioms \eqref{axiom.H}--\eqref{axiom.product}
above. Define $\mW$ to be the algebraic sum of the spaces
$\mW^{\mu}_{k} := q(S_{cl}^{\mu} ) + \mW^{-\infty}_{k}$, $\mu \in
\CC$, as above. Then $\mW$ is an algebra of (possibly) unbounded
operators on $\mH$ with dense common domain $\mW^{-\infty}\mH$ and
satisfying $\mW^{\mu}_{k} \mW^{\nu}_{l} \subset \mW^{\mu + \nu}_{k + l}$
Similarly, $\mW_{1,0}$, the union of the increasing sequence of
subspaces $\mW_{1,0}^{m} := q(S_{1, 0}^{m} ) + \mW^{-\infty}$, 
$m \in \ZZ$, is also an algebra and $\mW \subset \mW_{1, 0}$.
\end{proposition}

We shall endow the spaces $\mW^{\mu}_{k}$ and $\mW^m_{1, 0}$
introduced in the statement above with the quotient topology with respect
$q(S_{cl}^{\mu} ) \oplus \mW^{-\infty}_{k}\to \mW^{\mu}_{k}$
and $q(S_{1, 0}^{m}) \oplus \mW^{-\infty}\to \mW_{1, 0}^{m}$.
Let $\mW^{j} := \cup_{k \ge 1} \mW^{j}_{k}$. 
Similarly, we endow $\mW^{m}_{1,0,k} := q(S_{1,0}^{m} ) +
\mW^{-\infty}_{k}$ with the quotient topology. We clearly have
$\mW^{m}_{1,0} = \cup_{k} \mW^{m}_{1,0,k}$.

\begin{proof}\ 
First let us notice that Axiom \eqref{axiom.P} gives that
\begin{equation}\label{eq.prod.cl}
	q(S^{\mu}_{cl}) q(S^{\nu}_{cl}) \subset \mW^{\mu + \nu}_{1}.
\end{equation} 
Hence $\mW^{\mu}_{k}\mW^{\nu}_{l} \subset \mW^{\mu + \nu}_{k + l}$,
$\mu, \nu \in \CC \cup \{-\infty\}$, by Axiom \eqref{axiom.product'}
and the property $\mW^{-\infty}_{k} \mW^{-\infty}_{l} \subset
\mW^{-\infty}_{k + l}$ of an $LF$--algebra.
\end{proof}

\begin{definition}\label{def.extended}\ Assume the 
Axioms \eqref{axiom.H} -- \eqref{axiom.product}. Then
the algebra $\mW$ of Proposition \ref{prop.def.mW} above will be
called {\em an extended Weyl algebra}. The map $q : S^{\infty}_{1, 0}
\to \Hom(\mW^{-\infty}\mH, \mH)$ will be called the {\em quantization map}.
\end{definition}

Our definition differs from the definition in \cite{Guil85} in form
but not in substance.

\subsection{Conditions $(\sigma)$ and $(\psi)$}
For the construction of complex powers, we shall need our extended
Weyl algebra to satisfy two other conditions, Conditions $(\sigma)$
and $(\psi)$, which we now introduce.
\vspace*{2mm}
\\
$(\sigma)$ \ \vtop{%
{\em The map 
\begin{equation}\label{eq.sigma}
	S^{m}_{cl} \oplus \mW^{-\infty} \ni (a, T) \mapsto q(a) + T \in
	\mW^{m} := q(S^{m}_{cl}) + \mW^{-\infty} \subset \mW
\end{equation}
has kernel $\{ (a, -q(a)), a \in S^{-\infty}_{cl}\} \subset
S^{-\infty}_{cl}\oplus \mW^{-\infty}$.}  }
\vspace*{2mm}

If this is the case, the algebra $(\cup_{m \in \ZZ}
\mW^{m})/\mW^{-\infty}$ will be a topologically filtered algebra, in
the sense of \cite{BenameurNistor1, BenameurNistor2}.
\vspace*{2mm}

$(\psi)$\ \ {\em $\mW^{-\infty} = \mW^{-\infty}_{1}$} (so
$\mW^{-\infty}$ is actually a Fr\'echet algebra) {\em with topology
generated by a submultiplicative family of seminorms
$\|\;\cdot\;\|_{n}$, $n \ge 0$, such that for any operator of the form
$I + R$, $R \in \mW^{-\infty}$, that is invertible as a bounded
operator on $\mH$,
there is an $R_1\in \mW^{-\infty}$ with 
$(I + R)^{-1} = I + R_1$. }
\vspace*{2mm}

Recall that a {\em$\Psi^*$-algebra with unit} is a Fr\'echet algebra
with unit $\mB$ that is continuously embedded in $\mL(\mH)$, for some
Hilbert space $\mH$, such that $\mB^{-1} = \mL(\mH)^{-1} \cap \mB$.
(Here $\mB^{-1}$ denotes the group of invertible elements of $\mB$,
and the last condition means that $\mB$ is {\em spectrally invariant} in
$\mL(\mH)$.)  $\Psi^*$-algebras were introduced by Gramsch, see
\cite{gra84, lmn2, Schrohe}.

Thus, in Gramsch's terminology, condition $(\psi)$ is equivalent to
saying that $\mW^{-\infty} + \CC I$ is a $\Psi^*$-algebra whose
topology is generated by a multiplicative family of seminorms.  

It is useful now to recall an old result of Banach \cite{ban48},
which says that the inversion in a Fr\'echet algebra is continuous if,
and only if, its group of invertible elements is a
$G_{\delta}$-set. In particular, inversion is continuous on a
$\Psi^*$-algebra.

Although $(\psi)$ does not hold in some pseudo-differential algebras
that one constructs directly, such as the small b-calculus
$\Psi_{\mathrm{b}}^*(M)$, it does hold for their Guillemin completion
(see Proposition~\ref{prop.psi.star} and the example below).

\subsection{Examples}\ We now briefly discuss some examples. More
details will be provided in \cite{alnv2}.

\begin{example}\label{ex.sc}\ 
{\em The scattering calculus on $\RR^n$:}\ \cite{Cordes,
MelroseScattering, Parenti, ShubinDoc} Let $B$ be the closed unit ball
in $\RR^n$. The map $\Phi : \RR^d \ni x \mapsto x\<x\>^{-1} \in B$ defines
a diffeomorphism onto the interior of $B$ (the ``radial
compactification map'').

We shall take $M = B$, $V = B \times \RR^d$, $\mH = L^2(\RR^d)$, and
$q$ to be the Weyl quantization,
$$
	\big [q(a)u\big ](x) = (2\pi \imath)^{-d}\int_{\RR^d} \left (
	\int_{\RR^d} e^{\imath \<x-y, \xi\>} a(\Phi(\frac{x + y}{2}),
	\xi) u(y)\, dy \right ) d\xi.
$$
and 
$$
	\mW^{-\infty} = \mW^{-\infty}_{1} := q(S^{-\infty} ),
$$
with the induced topology. Let $\Delta = d^*d$ be the Laplace
operator, then we can take $R = e^{-t\Delta}$, $t > 0$. The results of
the papers mentioned above (see also \cite{HV, hor3, Schrohe,
Taylor2}) show that this is indeed an extended Weyl algebra.
\end{example}

Pseudodifferential operators on a compact, 
manifold without boundary form
an extended Weyl algebra as well.

\begin{example}\label{ex.comp}\ 
{\em Compact manifolds:}\ Let $M$ be a compact manifold.  We take $V =
T^*M$, $\mH = L^2(M)$, and let $\mW = \mW^{-\infty}_{1}$ be the
algebra of smoothing operators, $\mW=\Psi^{-\infty}(M) \simeq \CI(M
\times M)$ acting on $L^2(M)$ by convolution.  To define a 
suitable quantization map $q$, let us cover $M$ with finitely many coordinate
neighborhoods $U_\alpha \simeq \RR^d$. Let $q_{\alpha}$ be the
quantization map corresponding to this identification of $U_{\alpha}$
with $\RR^d$ (as in the example above) and let $\sum_{\alpha}
\phi_{\alpha}^{2}$ be a smooth partition of unity subordinated to
$U_{\alpha}$. Then we define
\begin{equation}\label{eq.quant.M}
	q(a) = \sum_{\alpha} \phi_{\alpha} q_{\alpha}(a)
	\phi_{\alpha}\,.
\end{equation}
Note that, this quantization map depends on the choice of the 
trivializations and of the cut-off functions.

We can again take $R = e^{-t\Delta}$, $t > 0$. It is a classical fact
that we obtain in this way an extended Weyl algebra (it also
follows, by localization, from the first example).
\end{example}

\begin{example}
Let $M$ be a smooth compact manifold with boundary and let
$x\in\CI(M)$ be a boundary defining function (so $x\geq 0$, $\partial
M$ is its zero set, and $dx\neq 0$ at $\partial M$). Let
$V={}^{\mathrm b}T^*M$, the dual bundle of ${}^{\mathrm{b}}TM$, whose
smooth sections are the smooth vector fields on $M$ that are tangent
to $\partial M$. The small calculus $\mW^m=\Psi_{\mathrm{b}}^m(M)$ is
defined in \cite{MelroseAS} as the algebra of operators whose Schwartz
kernel is conormal (of order~$m$) to the lifted diagonal on the blown
up space $M^2_{\mathrm{b}}=[M\times M;\partial M\times\partial M]$ and
vanishes to infinite order at the left and right faces, i.e.\ at the
lifts of $\partial M\times M$ and $M\times\partial M$ to
$M^2_{\mathrm{b}}$.  The elements of $\mW^{-\infty}$ are those
operators whose Schwartz kernel is in addition smooth. (It is easy to
write down a quantization map $q$ as in the previous example.)
Moreover, let $\mH=L^2_{\mathrm{b}}(M)=x^{1/2} L^2(M)$, be the natural
$L^2$-space with respect to b-densities, defined using $V$.  Then $M$
does not satisfy~$(\psi)$: typically the inverses of operators in
$I+\mW^{-\infty}$, even if they exist in $\mL(\mH)$, have Schwartz
kernels that vanish only to finite order at the left and right
faces. However, it does satisfy all other axioms, including
$(\sigma)$. Then Proposition~\ref{prop.psi.star} will show that the
Guillemin completion $\overline{\mW}^{-\infty}$, i.e.\ the completion
of $\mW^{-\infty}$ in the topology of the operator norms
$H_{\mathrm{b}}^s(M)\to H_{\mathrm{b}}^{s'}(M)$, satisfies all
axioms. This illustrates the power, and limitations, of our results:
for a b-metric $g$, we show that $(\Delta_g+1)^z$ is in
$\overline{\mW}^{2z}$.  However, it certainly does not lie in
$\mW^{2z}$, although it does lie in a `large b-calculus' \cite{Loya1}. 
Note that in a similar closures of the $b$-calculus have been studied 
in detail by Mantlik \cite{man95}.
\end{example}

Manifolds with bounded geometry can also be included in our framework.

\begin{example}\label{ex.bdg}\
Let $M$ be a manifold with bounded geometry and injectivity radius $r$.
By definition of ``bounded geometry'' $r>0$.  Take $\mW^{\mu} =
B\Psi^{\mu}(M)$ \cite{kordyukov, Shubin}, and $\mW^{\mu}_{k}$ be the
operators with support in the set $M^2_{kr/2}$ of pairs of points
$(x,y) \in M^2$ at distance $\le kr/2$.  The topology on
$\mW^{-\infty}_{k}$ is as defined in \cite{kordyukov, Shubin}, which
in fact coincides with the restriction of the Guillemin topology (see
Proposition \ref{prop.psi.star} below) on this subspace.
\end{example}

Our main interest is in manifolds with a Lie structure at infinity,
which generalize the first two examples above and refine the algebras
$B\Psi^\infty$ of the last example.

\begin{example}\label{ex.Lie}\ 
Let $(M,A)$ ba a manifold with a Lie structure at infinity and
$\mW^{\mu} = \Psi^{\mu}(M,A)$ be the algebras of pseudodifferential
operators canonically associated to it, see \cite{aln1, alnv2}. Let
$r$ be the injectivity radius of $M$. We let as above $\mW^{\mu}_{k}$ be
the operators with support in the set $M^2_{kr/2}$ introduced in the
previous example. This example will be discussed in detail in
\cite{alnv2}.
\end{example}

Yet another example is obtained by considering the so called
``$\Phi$--calculus'' of Mazzeo and Melrose, see
\cite{MaMe}.

\section{First consequences of the definitions\label{Sec.FC}}

We now establish the first consequences of the axioms. One of our
goals is to construct Sobolev spaces, in the next section.
Only Axioms \eqref{axiom.H}--\eqref{axiom.cont}
are needed to define Sobolev spaces (Axiom~\eqref{axiom.product} will
not be needed, so in this section and next two sections we do not assume it).  

As above, let $\mW = \sum \mW^{\mu}_{k}$, $\mW^{\mu}_{k} := q(S_{cl}^{\mu} ) +
\mW^{-\infty}_{k}$, $\mu \in \CC$, and $\mW_{1,0} = \sum \mW^m_{1,0}$, 
$\mW^m_{1,0} := q(S^m_{1,0}) +
\mW^{-\infty}$, $m\in \ZZ$.

As in \cite{Guil85}, an element $A \in M_N(\mW^\mu)$ is called {\em
elliptic}, if there exists $B \in \mW^{-\mu}$ such that $AB - I \in
M_N(\mW^{-1})$.  Note that our assumptions  (especially those for $\cP_0$) give
$$
	AB - I \in M_N(\mW^{-1}) \Leftrightarrow BA - I \in
	M_N(\mW^{-1}).
$$
Similarly, a symbol $a \in M_N(S^{m}_{1, 0})$ is called {\em elliptic}
if there exists $b \in M_N(S^{-m}_{1, 0})$ such that $ab - 1 \in
S^{-1}_{1, 0}$, as usual.

\begin{lemma}\label{lemma.sh}\ 
Let $T \in \mW^{m}_{1, 0}$ and $\sh{T} := T^*\vert_{\mW^{-\infty}
\mH}$. Then $\sh{T} \in \mW^{\overline{m}}_{1, 0}$ and $T \mapsto \sh{T}$
defines a conjugate-linear 
involution on $\mW_{1, 0}$ satisfying $\sh{(\mW^{\mu}_{k})} \subset
\mW^{\overline{\mu}}_{k}$. Moreover, $q(\overline{a}) = \sh{q(a)}$ and
is hence contained in $q(a)^*$, the adjoint of $q(a)$, for any $a \in
S^{m}_{1, 0} = S^{m}_{1, 0}(V)$.
\end{lemma}

\begin{proof}\ 
Let $a = b + \imath c \in S^{m}_{1, 0} = S^{m}_{1, 0}(V)$ with $b$ and
$c$ real valued. Then
\begin{equation}\label{eq.adj}
	(q(a) \xi, \eta) = (q(b)\xi, \eta) + \imath (q(c)\xi, \eta) =
	(\xi, q(b)\eta) + \imath (\xi, q(c)\eta) = (\xi,
	q(\overline{a})\eta),
\end{equation}
for any $\xi, \eta \in \mW^{-\infty}\mH$, since $q(b)$ and $q(c)$ are
symmetric by the Axiom~\eqref{axiom.q}. Let $P \in \mW_{1, 0}$, then
$P = q(a) + T$, with $a \in S^{m}_{1, 0}$ and $P \in \mW^{-\infty}$.
Then Equation \eqref{eq.adj} implies that $P^*\xi = (q(\overline{a}) +
T^*) \xi$, for any $\xi \in \mW^{-\infty}\mH$. Thus $\sh{P} \in
\mW_{1, 0}$ and $\sh{(\mW^{\mu})} = \mW^{\overline{\mu}}$. Moreover,
$\sh{(\sh{P})} = q(a) + T = P$ and hence $\sh{}$ is an involution of
$\mW_{1, 0}$ and $\mW$. Also,
$$
	(\sh{(P_1 P_2)}\xi, \eta) = ( (P_1P_2)^*\xi, \eta) = (\xi ,
	P_1P_2\eta) = (\sh{P_2}\sh{P_1}\xi, \eta),
$$
for any $\xi, \eta \in \mW^{-\infty}\mH$, and hence $\sh{}$ is
anti-linear.  It is easy to check that $\sh{\lambda T} =
\overline{\lambda} \sh{T}$, and hence $\sh{}$ is conjugate-linear as
well.
\end{proof}

This lemma then gives the following.

\begin{proposition}\label{prop.cont}\ 
(i) There exists a continuous seminorm 
$p'$ on $S^{0}_{1, 0} := S^{0}_{1, 0}(V)$ such that $\| q(a) \| \le
p'(a)$, the norm here being the norm of bounded operators on $\mH$.

\noindent (ii) The map $S^{0}_{1,0} \ni a \mapsto q(a) \in \mL(\mH)$ is well defined
and continuous.

\noindent (iii) If $a \in 
S^m_{1,0}$, $m \ge 0$, is elliptic and real valued,
then $q(a)$ is essentially self-adjoint.
\end{proposition}

\begin{proof}\ 
For proving (i) and (ii) we will use  the symbolic calculus and H\"ormander's trick 
\cite{fio, hor3}

For $a \in S^{0}_{1, 0}$ we take a constant $M\ge |a| + 1$, $M \in \RR$. 
We can assume $a$ real, for simplicity. 

Then we define
\begin{equation}\label{eq.b}
        b_0 = (M^2 - a^2)^{1/2} \in S^{0}_{1, 0}
\end{equation}
Let $r_0 = M^2 - \mP(a,a) - \mP(b_0,b_0)$. Then $r_0 \in
S^{-1}_{1,0}$, $q(r_0) - (M^2 - q(a)^2 - q(b_0)^2) \in \mW^{-\infty}$
and $r_0$ will depend continuously on $a$ because $\mP(a',b') - a'b'
\in S^{ - 1}_{1,0}$, if $a' \in S^{0}_{1,0}$ and $b' \in
S^{0}_{1,0}$, and the induced map is continuous
(Axiom~\eqref{axiom.P}).

Let $b_1 = b_0 + b'_1$, with
$$
	b'_1 = b_0^{-1} \big  [r_0 - M^2 - \mP(a,a) \big ]/2 \in S^{-1}_{1,0}
$$ 
and $r_1 = M^2 - \mP(a,a) - \mP(b_1,b_1)$. Then $r_1 \in
S^{-2}_{1,0}$, $q(r_1) - (M^2 - q(a)^2 - q(b_1)^2) \in \mW^{-\infty}$
and $r_1$ will depend continuously on $a$. Iterating this
construction, we see that we can assume that we have constructed $b_d
\in S^0_{1,0}$, $r_d \in S^{-d-1}_{1,0}$, both depending continuously
on $a$, such that $r_d = 0$ if $a=0$ and
$$
	R := q(r_d) - (M^2 - q(a)^2 - q(b_d)^2) \in \mW^{-\infty}
$$ 
also depends continuously on $a$.

Now we use Axiom~\eqref{axiom.cont} to conclude that for sufficiently large
$d$ the operator $q(a)$ is bounded and $\|q(a)\|$ is continuous in $a$.

We then obtain
\begin{equation}\label{eq.middle}  
 	\|q(a)\xi\|^2 = (q(a)\xi, q(a)\xi) = M^2\|\xi\|^2 + (R\xi, \xi)
	-(q(r_d)\xi, \xi) - \|q(b)\xi\|^2,      
\end{equation}
for any $\xi \in \mW^{-\infty}\mH$, which is dense in $\mH$, by
Axiom~\eqref{axiom.R}. This gives
$$
	\|q(a)\|^2 \le (\|a\|_{\infty} + 1)^2 + \|R\| + \|q(r_d)\| =
	C(a),
$$ 
where $C(a)$ is a continuous function of $a \in S^0_{1,0}$, by our
discussion above. This proves (ii). The first part, (i), is then a
direct consequence of (ii).

To prove (iii), we can assume $m > 0$, otherwise the statement is
trivial, because $q(a)$ is bounded in that case. Since $q(a)$ is
symmetric by the Axiom \eqref{axiom.q}, it is enough to check that
$(q(a) \pm \imath I)(\mW^{-\infty}\mH)$ are dense 
\cite[Theorem VIII.3 and Corollary]{ReedSimon}. Using the asymptotic
completeness of the space of symbols, we can find $b \in
S^{-m}_{1,0} $ such that $ab - 1 \in S^{-1}_{1,0} $ and $(q(a) + \imath
I) q(b) - I \in \mW^{-\infty}$.

Choose $R = R^* \in \mW^{-\infty}$ to
be an injective map. Then
$$
	P = (q(a) - \imath I) (R^*R + q(\overline{b})q(b)),
$$
also satisfies $(q(a) + \imath I) P = I + T,$ for some $T \in
\mW^{-\infty}$. Note that
$$
	(R^*R + q(\overline{b})q(b)\xi, \xi) \ge (R\xi, \xi) > 0
$$
if $\xi \in \mW^{-\infty}\mH$, $\xi \not = 0$, by Lemma
\ref{lemma.sh}.  Now, if $(q(a) + \imath I)(\mW^{-\infty})\mH$ is not
dense, we can find $\eta \in \mH$, $\eta \not = 0$, perpendicular to
$(q(a) + \imath I)P (\mW^{-\infty}\mH)$.

This gives that $(I + T)^*\eta$ is perpendicular to
$\mW^{-\infty}\mH$.  Since the latter 
is dense in $\mH$ by Axiom
\eqref{axiom.R}, we obtain that $(I + T)^*\eta = 0$, and hence $\eta =
- T^*\eta \in \mW^{-\infty}\mH$ is in the domain of all operators in
$\mW$. But then
$$
	(\eta, (q(a) + \imath I)P \zeta ) = ( P^*(q(a) - \imath
	I)\eta, \zeta ) = 0
$$
for all $\zeta \in \mW^{-\infty}\mH$, which implies as above that
$$
	(R^*R + q(b)^{*}q(b)) (q(a) + \imath I) (q(a) - \imath I)\eta
	= 0.
$$
Now this is a contradiction, because all operators in the above
product are injective. (Note that we have used above $(q(a) + \imath
I)^*\eta = (q(a) - \imath I)\eta$ and similarly for the other
operators above. The main point of the proof was to show that this is
indeed possible.)

The proof for the other choice of sign in $(q(a) \pm \imath I)$ is the
same, and hence our proof of (iii) is complete.
\end{proof}

For further reference, let us mention here a consequence of being elliptic and 
invertible in the algebra $\mW^{\infty}$.

\begin{proposition}\label{prop.invert}\ 
Let $m \ge 0$, $A \in M_N(\mW^{m})$, and $B \in M_N(\mW^{-m})$ be such
that $AB = BA = I$ in $M_N(\mW)$. Then $\overline{A}$, the closure of
$A$, is invertible with inverse $\overline{B}$. In particular, the
domain of $\overline{A}$ is $\overline{B}\mH$.
\end{proposition}

\begin{proof}\ 
Use the boundedness of $B$ and the definitions.
\end{proof}

\section{Sobolev spaces\label{Sec.AS}}

We now define some abstract Sobolev space $H^{(s)}$ as in Guillemin's
paper.  Recall from the previous section that only the Axioms
\eqref{axiom.H}--\eqref{axiom.cont}
are assumed in this and next section.

As in Guillemin's paper, we define the Sobolev space $H^{(s)}$, $s \ge
0$, to be the domain of (the closure of) $P$, where $P \in \mW^{s}$ is
elliptic. This definition is independent of the choice of $P$, by
Axioms \eqref{axiom.P} and \eqref{axiom.cont}, because if $P_1$ is
another such elliptic operator, then we can find $Q \in \mW^0$ and $R
\in \mW^{-\infty}$ such that $P_1 = QP + R$ and both $Q$ and $R$ are
bounded.  Let $r \in S^{1}_{cl}$ be 
a symbol with 
\begin{equation}\label{def.glob.r}
  r(\xi)= \|\xi\|\qquad \mbox{for all $\xi\in V$ with $\|\xi\|\geq 2$}
\end{equation}
and 
\begin{equation}\label{eq.Ps}
	P_s := \frac{1}{2} \big (I + q(r^{s/2})^{2}) \big ).
\end{equation}
We endow $H^{(s)}$ with the norm defined by
$$
	\|f\|_{(s)}^2 := \| P_s f\|^2,
$$
$f \in \mH$ in the domain of (the closure of) $P_s$.  If $s < 0$, we
define $P_{s} := P_{-s}^{-1}$ and $H^{(s)}$ to be the completion of
$\mH$ in the norm $\xi \to \|\xi\|_{(s)} := \|P_{s}\xi\|$. Then the
inner product of $\mH$ extends to a bilinear pairing $H^{(s)} \otimes
H^{(-s)} \to \CC$ that identifies $H^{(-s)}$ with the dual of
$H^{(s)}$ as in \cite{LauterMoroianu}, for example.

\begin{proposition}\label{prop.bd.Sob}\ 
(i)\ An operator $T \in \mW$ defines a continuous map $T : H^{(s)} \to
H^{(s')}$ if, and only if, $P_{s'} T P_{-s}$ extends to a bounded
operator on $\mH$.

(ii)\ In particular, any $P \in \mW^z$ extends by continuity to a
bounded map $H^{(s)} \to H^{(s-r)}$ for $r = \re(z)$.
\end{proposition}

\begin{proof}\ 
We have from definition that $P_s : H^{(s)} \to \mH$ is an isometric
isomorphism and that $P_{-s} = (P_{s})^{-1}$ for all $s$. (We first
replace all operators by their closures, note however that no
confusion can arise, in view of Proposition \ref{prop.invert}.) This
takes care of (i).

Let $z = r + i\mu$, $r, \mu \in \RR$. The second part is proved directly using
the symbolic calculus for the four possible signs of $s$ and
$s':= s-r$.

If $-s, s' \ge 0$, then $P_{s'}TP_{-s} \in \mW^{\imath \mu}$ is a
bounded operator.

If $s, s' \ge 0$, the symbolic calculus tells us that we can write $T
= T_1 P_s + R$, where $T_1 \in \mW^{-s'}$ and $R \in
\mW^{-\infty}$. Then
$$
	P_{s'} T P_{-s} = P_{s'} T_1 P_s P_{-s} + P_{s'} R P_{-s} =
	P_{s'} T_1 + P_{s'} R P_{-s}
$$
is bounded because $P_{s'}T_1 \in \mW^{\imath \mu}$ and $P_{s'}RP_{-s}
\in \mW^{-\infty}$ are bounded.

If $-s, -s' \ge 0$, we proceed similarly by writing $T = P_{-s'}T_1 +
R$ with $T_1 \in \mW^{s + \imath \mu}$ and $R \in \mW^{-\infty}$.

If $s, -s' \ge 0$, we proceed similarly by writing $T = P_{-s'}T_1
P_s + R$ with $T_1 \in \mW^{\imath \mu}$ and $R \in \mW^{-\infty}$.
\end{proof}

\section{Holomorphic families of operators\label{Sec.HFO}}

Recall from the previous sections that only the Axioms
\eqref{axiom.H}--\eqref{axiom.cont} are assumed in this section.

Let $\mO(\Omega)$ be the space of holomorphic functions on some open
subset $\Omega \subset \CC$. If $X$ is a Fr\'echet space, we denote by
$\mO(\Omega, X)$ the space of holomorphic (or complex differentiable)
functions on $\Omega$ with values in $X$.

{\em We shall assume from now on throughout the paper that our algebra
$\mW$ satisfies condition $(\sigma)$ of Section \ref{Sec.EWA}.}

Condition $(\sigma)$ allows us to define the principal symbol
$\sigma^{(s)}(T)$ of an operator $T \in \mW^{s}$ as in the classical
case of algebras of pseudodifferential operators. Namely,
$$
	\sigma^{(s)} : \mW^s \to S_{cl}^{s}/S_{cl}^{s-1}
	\subset \CI(V \smallsetminus 0) 
$$ 
satisfies $\sigma^{(s)}(q(a)) = a + S_{cl}^{s-1}(V)$, which determines
it completely. (Let $V \smallsetminus 0$ be the set of non-zero
vectors of the vector bundle $V$. We identify
$S_{cl}^{s}/S_{cl}^{s-1}$ with the space of smooth, homogeneous of
order $s$ functions on $V \smallsetminus 0$.) The maps $\sigma^{(s)}$
and $q$ immediately extend to $M_N(\mW^s)$. Then $A$ is elliptic if,
and only if, $\sigma^{(s)}(A)$ is invertible everywhere on 
$V\setminus\{0\}$.

We now establish an analogue of the asymptotic completeness for
holomorphic families of operators in $\mW$. We shall write $M_N(X)$ for
$M_N(\CC) \otimes X$, for any Fr\'echet space $X$.

\begin{proposition}\label{prop.h.a.c}\
(i)\ Let $a_i \in \mO(\Omega, M_N(S_{cl}^{\mu - i} ))$, $i = 0, 1, 2,
\ldots$, be holomorphic functions on some open subset $\Omega \in
\CC$. Then there exists a holomorphic function $a \in \mO(\Omega,
M_N(S^{\mu}_{cl}))$ such that
for all $k=1,2,\dots$ 
   $$a(z) - \sum_{i = 0}^{k - 1} a_i(z) \in M_N(S_{cl}^{\mu - k} ).$$

(ii)\ For any $F \in \mO(\Omega, M_N(\mW^\mu_{k}))$ there exists $f \in
\mO( \Omega, M_N(S_{cl}^\mu))$ such that $F(z) - q(f(z)) \in
M_N(\mW^{-\infty}_{k})$ is a holomorphic function.

(iii)\ Similarly, let $A_i \in \mO(\Omega, M_N(\mW^{\mu - i}_{k}))$, $i =
0, 1, 2, \ldots$, be holomorphic functions on some open subset $\Omega
\in \CC$. Then there exists a holomorphic function $A \in
M_N(\mO(\Omega,\mW^{\mu}_{k}))$ such that 
for all $l=1,2,\dots$ 
  $$A(z) - \sum_{i = 0}^{l - 1} A_i(z)\in M_N(\mW^{\mu - l}_{k}).$$
\end{proposition}

If $a$ is as in (i) above, we shall say that $a$ is an {\em asymptotic
sum} of the sequence $a_i$, $i \ge 0$. Similarly, if $A$ is as in
(iii) above, we shall say that $A$ is an {\em asymptotic sum} of the
sequence $A_i$, $i \ge 0$.

\begin{proof}\ 
Assume $\mu = 0$, $N = 1$, and $\Omega = \CC$, for simplicity.

Let $B$ be the set of vectors of length at most one in $V$, which we
identify with the radial compactification of $V$ using the map $V \ni
\xi \longmapsto \<\xi\>^{-1}\xi \in V$.  This identifies $S_{cl}^0 $
with $\CI(B)$ and we let
$$
	p_n(a) = {\max}_{l=0}^n \|\nabla^la\|_\infty
$$ 
be the standard seminorms defining the topology on $\CI(B^*)$.

We observe that $S_{cl}^{-n-1} $ is contained in the closure of
$S^{-\infty} $ in the topology defined by the seminorm $p_n$. By
applying this to the Taylor coefficients of $a_j(z)$, for any $j$, we
obtain that there exist holomorphic functions $r_j \in
\mO(\Omega, S^{-\infty})$ such that $p_j(a_{j+1}(z) - r_{j+1}(z)) \le
2^{-j}$ for any $|z| \le j$. This proves that the series $\sum_{j}
(a_{j}(z) - r_{j}(z))$ is uniformly convergent on the compact subsets
of $\CC$. Let $a(z)$ be its sum. Then $a \in \mO(\Omega, S_{cl}^0 )$.
This proves (i).

To prove (ii), proceed as follows. Let $\chi$ be a smooth function on
$V$, $\chi = 1$ outside the unit ball and $\chi = 0$ on the set of
vectors of length $\le 1/2$. Define $a_0(z) = \chi
\sigma^{(0)}(A(z))$, which is a holomorphic function with values in
$S^{0}_{cl} $.  We define $a_{n}(z)$, $n \ge 1$, by induction by
$$
	a_{n}(z) = q \Big[ \chi \sigma^{(-n)}\big( A(z) - \sum_{j =
	0}^{n-1}q(a_j(z)) \big ) \Big]\,.
$$
Then use (i) to construct $a \in \mO(\Omega, S_{cl}^0 )$ such that
$a(z) - \sum_{i = 0}^{k - 1} a_i(z) \in S_{cl}^{- k} $. Then $a$
satisfies the requirements of (ii).

We now turn to (iii). First, using (ii), we find $a_i \in \mO(\Omega,
S^{-i}_{cl} )$ such that $A_i(z) - q(a_i(z)) \in \mW^{-\infty}$. Then
we use (i) for the sequence $a_i$ to find $a \in \mO(\Omega,
S^{0}_{cl} )$ such that $a(z) - \sum_{i = 0}^{k - 1} a_i(z) \in
S_{cl}^{-k}.$ We can let then $A = q(a)$.
\end{proof}

We shall write
$$
	\mO(\Omega, \mW^{\mu}) := \cup \mO(\Omega, \mW^{\mu}_{k}).
$$
Standard reasonings then give the following result.

\begin{corollary}\label{cor.elliptic}\ 
(i) Let $A \in \mO(\Omega, \mW^{\mu})$ be a holomorphic, pointwise
elliptic family. Then there exists $B \in \mO(\Omega, \mW^{-\mu})$
such that $AB - I, BA - I \in \mO(\Omega, \mW^{-\infty})$.

(ii) There exists $b \in \mO(\CC, S_{cl}^0)$, $b(0) = 1$, such that
$q(r^z) q(r^{-z}b(z)) - I \in \mW^{-\infty}$ and, similarly,
$q(r^{-z}b(z)) q(r^z) - I \in \mW^{-\infty}$ with $r$ defined via ~\eqref{def.glob.r}.
\end{corollary}

\begin{proof}\ 
(i) follows right away from Proposition \ref{prop.h.a.c} using
standard arguments \cite{hor3, Taylor2}. Namely, let $a \in
\mO(\Omega, S^{\mu}_{cl} )$ be such that $A(z) - q(a(z)) \in
\mW^{-\infty}$ for all $z \in \Omega$. Choose $\chi \in \CI(V)$ to be
equal to zero in a neighborhood of $0$ containing all zeroes of
$a(z)$, for all $z$, and equal to $1$ outside some larger neighborhood
of $0 \in V$.  Let $b(z) := \chi a(z)^{-1}$, which is defined to be
zero where $a(z) = 0$. Let $B_0(z) = q(b(z)) \in \mO(\Omega,
\mW^{-\mu})$. Then $AB_0 = I + R_{1}(z)$, where $R_1 \in \mO(\Omega,
\mW^{-1})$. We can find $R_{1n} \in \mO(\Omega, \mW^{-n}_{1})$ such
that $R_{1n}(z) - R_1(z)^n \in \mW^{-\infty}$.  Define $B$ then to be
an asymptotic sum of the sequence $(-1)^{n} B_{0} R_{1n}$ (see
Proposition \ref{prop.h.a.c} (iii)).

We now apply (i) to the holomorphic, elliptic family $A(z) = q(r^{z})
q(r^{-z})$ to deduce that there exists $a \in \mO(\Omega, S_{cl}^0 )$
such that $q(r^{z})q(r^{-z})q(a(z)) - I \in \mW^{-\infty}$. Then
$q(r^{-z})q(a(z)) = q(r^{-z}b(z)) + r(z)$, where $b \in \mO(\Omega,
S_{cl}^0 )$ and $r \in \mO(\Omega, \mW^{-\infty})$, by Axiom
\eqref{axiom.P}.
\end{proof}

\section{Spectral invariance\label{Sec.SI}}

In this section we prove the continuity of the multiplication
$\mW^\mu_k \times \mW^\nu_l \to \mW^{\mu + \nu}_{k + l}$, we
introduce the Guillemin completion of an extended Weyl algebra,
which, in fact, is again an extended Weyl algebra satisfying
also Condition $(\psi)$. Then we establish the spectral invariance
of extended Weyl algebras satisfying conditions $(\sigma)$ and $(\psi)$.

{\em From now on and throughout the rest of the paper, $\mW$ will
denote an extended Weyl algebra, i.e. we assume
Axioms~\eqref{axiom.H}--\eqref{axiom.product}. }

As before let $r\in S_{cl}^{1}$ be a symbol satisfying
\eqref{def.glob.r}.  We now establish the partial continuity of the
multiplication on an extended Weyl algebra.

\begin{proposition}\label{prop.cont.m}\
The map $\mP : S^{\mu}_{cl} \times S^{\nu}_{cl} \to S^{\mu +
\nu}_{cl}$ and the multiplication maps
$$
	\mW^\mu_{k} \times \mW^{-\infty}_{l} \to \mW^{-\infty}_{k +
	l},\quad \mW^{-\infty}_{k} \times \mW^{\mu}_{l} \to
	\mW^{-\infty}_{k + l},\quad \mW^\mu_{k} \times \mW^{\nu}_{l}
	\to \mW^{\mu + \nu}_{k + l},
$$ 
where $\mu, \nu \in \CC$ and $k, l = 1, 2, \ldots$, are well
defined and continuous.
\end{proposition}

\begin{proof}\ 
Axiom \eqref{axiom.H} tells us that $\mW^{-\infty} = \cup_{k \ge 1}
\mW^{-\infty}_{k}$ is an $LF$--algebra, and hence the multiplication
map
$$
	\mW^{-\infty}_{k} \times \mW^{-\infty}_{l} \to
	\mW^{-\infty}_{k + l}
$$ 
is continuous. Axiom \eqref{axiom.product} shows that the
multiplication map
$$
	S^{\mu}_{cl} \times \mW^{-\infty}_{k} \ni (a, T) \to q(a)T \in
	\mW^{-\infty}_{k + 1}
$$ 
is continuous. Hence the induced map
$$
	(S^{\mu}_{cl} \oplus \mW^{-\infty}_{k} ) \times
	\mW^{-\infty}_{l} \to \mW^{-\infty}_{k + l}
$$
is also continuous (note that we used $k + l \ge k + 1$).  The
continuity of the first two maps is then a consequence of the
definition of $\mW^{\mu}_{k} := q(S^{\mu}_{cl}) +
\mW^{-\infty}_{k}$. For the continuity of the right multiplication, we
also use the continuity of $T \mapsto T^*$ on $\mW^{-\infty}_{k}$, see
Axiom \eqref{axiom.H}.

Let $\mR(a, b) = q(a) q(b) - q(\mP(a,b)) \in \mW^{-\infty}_1$.  Axiom
\eqref{axiom.P} and the definition of the topology on $\mW^{\mu +
\nu}_{1} := q(S^{\mu + \nu}_{cl}) + \mW^{-\infty}_{1}$ shows that the
maps
\begin{multline*}
\begin{gathered}
	S^{\mu}_{cl} \times S^{\nu}_{cl} \ni (a,b) \to \big( \mP(a,b),
	\mR(a,b) \big) \in S^{\mu + \nu}_{cl} \oplus \mW^{-\infty}_{1}
	\quad \text{and}\\ S^{\mu + \nu}_{cl} \oplus \mW^{-\infty}_{1}
	\ni (a, T) \to q(a) + T \in \mW^{\mu + \nu}_{1}
\end{gathered}
\end{multline*}
are well defined and continuous. Their composition is the product
$q(a)q(b)$, which hence depends continuously on $a$ and $b$ in the
above spaces of symbols. Then the same arguments as above, together
with the continuities already proved, show that the multiplication
$$
	\mW^\mu_{k} \times \mW^{\nu}_{l} \to \mW^{\mu + \nu}_{k + l}
$$ 
is well defined and continuous.
\end{proof}

Denote by $\|\;\cdot\;\|_n$ the operator norms $H^{(-n)} \to H^{(n)}$.
Let $\overline{\mW}^{-\infty}$ be the closure of $\mW^{-\infty}$ in
the topology defined by the family of norms $(\|\;\cdot\;\|_n)$, $n
\in \ZZ_+$.  Also, let
$$
	\overline{\mW} :=  \mW +  \overline{\mW}^{-\infty}.
$$
We shall call $\overline{\mW}$ and $\overline{\mW}^{-\infty}$ the {\em
Guillemin completions} of $\mW$ and $\mW^{-\infty}$, respectively. We
shall also write $\overline{\mW}^s := \mW^s + \overline{\mW}^{-\infty}$.

\begin{proposition}\label{prop.psi.star}\ 
Let $\mW$ be an extended Weyl algebra.  Then the Guillemin completion
$\overline{\mW}$ is an extended Weyl algebra that satisfies $(\psi)$
with respect to the family $(\|\;\cdot\;\|_n)$, ${n \in \ZZ_+}$, of
operator norms $H^{(-n)} \to H^{(n)}$.
\end{proposition}

(It is not assumed that $\mW$ satisfies $(\psi)$, nor that it is a
Fr\'echet algebra. The conclusion is that $\overline{\mW}$ has these
properties. Also, in applications, one will have to check that
$\overline{\mW}$ satisfies Condition $(\sigma)$ in order to be able to
construct complex powers. However, this condition is automatically
satisfied in practice for algebras of pseudodifferential operators.)

\begin{proof}\  
We shall check that $\overline{\mW}$ satisfies all the axioms defining
an extended Weyl algebra. To start with, the axioms \eqref{axiom.H},
\eqref{axiom.R}, and \eqref{axiom.cont} are automatically satisfied.

Recall that $\mL(X,Y)$ denotes the space of continuous linear maps
between two locally convex vector spaces $X$ and $Y$. If $X$ and $Y$
are Banach spaces, we denote by $\|\;\cdot\;\|_{\mL(X,Y)}$ the
canonical norm on $\mL(X,Y)$. In particular,
$$
	\|\; \cdot \;\|_n = \|\; \cdot \;\|_{\mL(H^{(-n)}, H^{(n)})}
$$ 
and $\| T \|_{n} = \| P_{n} T P_{n}\|_{0}$, where $P_{s}:= \frac{1}{2}
\Big (I + q(r^{s/2})^{2}) \Big )$ is as introduced in Equation
\eqref{eq.Ps}. In particular, all semi-norms $\|\;\cdot\;\|_{n}$ are
continuous on $\mW^{-\infty}_{k}$, and hence the inclusion map
\begin{equation}
	\mW^{-\infty}_{k} \hookrightarrow \overline{\mW}^{-\infty}
\end{equation}
is continuous, for any $k$. The continuity of this map is enough to
conclude that Axioms \eqref{axiom.q}, \eqref{axiom.P}, and
\eqref{axiom.product} are satisfied.

Let us notice that $\| T \|_{\mL(H^{(s)}, H^{(s')})} \le \| T \|_{n}$,
whenever $s \ge -n$ and $s' \le n$. This gives
$$
	\|T_1 T_2\|_{n} \le \| T_1 \|_{\mL(H^{(0)}, H^{(n)})} \| T_2
	\|_{\mL(H^{(-n)}, H^{(0)})} \le \| T_1 \|_{n} \| T_2 \|_{n}.
$$
Thus the family of seminorms $\|\; \cdot \;\|_n$ is submultiplicative.

The only thing that is left to prove is that the algebra 
$\CC + \overline{\mW}^{-\infty}$ is spectrally invariant in $\mL(\mH)$
(that is, that $\CC + \overline{\mW}^{-\infty} \cap \mL(\mH)^{-1} = 
(\CC + \overline{\mW}^{-\infty})^{-1}$), but this follows from Proposition 
\ref{semi}
due to the definition of the Sobolev spaces..
\end{proof}

It is worth pointing out that in a different axiomatic setting spectral invariance
of the closure of an algebra of bounded operators with respect to the family 
$\|\; \cdot \;\|_n$ of seminorms has been studied before by Mantlik \cite{man99}.

We shall need the following lemma.

\begin{lemma}\label{lemma.inv.o}\ \ 
If $\mW$ is an extended Weyl algebra satisfying $(\psi)$ and $w \in
\CC$, then we can find $Q \in \mW^{w}$ such that $Q^{-1} \in
\mW^{-w}$.
\end{lemma}

\begin{proof}\ 
Let $b \in \mO(\CC, S_{cl}^0)$, $b(0) = 0$, such that $R(z) :=
q(r^{z}) q(r^{-z}b(z)) - I \in \mW^{-\infty}$, as in Corollary
\ref{cor.elliptic}(ii). Then $R \in \mO(\CC, \mW^{-\infty})$ is
holomorphic and $R(0) = 0$. Hence $I + R(z) = q(r^{z}) q(r^{-z}b(z))$
is invertible on $\mH$ for $|z| < \epsilon$, for $\epsilon$
small. Then $q(r^{z})^{-1} = q(r^{-z}b(z)) (I + R(z))^{-1} \in \mW^{-z}$, 
for $|z| < \epsilon$ small.

Chose now $k \in \NN$ large enough so that $|w| < k\epsilon$. 
Then we can take $Q = q(r^{w/k})^k$.
\end{proof}

We then obtain

\begin{theorem}\label{theorem.invert}\ \ 
Let $\mW$ be an extended Weyl algebra satisfying $(\psi)$. Also, let
$P \in M_N(\mW^z)$, $\re(z) \ge 0$. We replace $P$ with its closure,
if necessary. Assume that $P$ is elliptic and invertible as a
(possibly) unbounded operator on $\mH$, then $P^{-1}$
is (the closure of an element) in $M_N(\mW^{-z})$.
\end{theorem}

\begin{proof}\ Below, we shall replace all operators in
$\mW$ by their closures, when necessary. We shall assume $N = 1$, for
simplicity.

Write $z = a + b\imath$, $a, b \in \RR$,
and let $Q \in \mW^{-z}$ be such that $Q^{-1} \in \mW^{z}$,
which is possible by Lemma \ref{lemma.inv.o}.  Then $P_1 := PQ \in
\mW^{0}$ is elliptic, injective, and bounded. The ellipticity of
$Q^{-1}$ guarantees the range of (the closure of) $Q$ is $H^{(a)}$
and the ellipticity of $P$ guarantees that the domain of $P$ is also
$H^{(a)}$. Thus $P_1 = PQ : \mH \to \mH$ is also surjective, and hence
it is invertible (on $\mH = H^{(0)}$).  We can hence replace $P$ by
$P_1$ and hence assume that it is bounded.

Let then $P \in \mW^{0}$ be invertible. At this point the proof can be
completed using \cite{lmn2}, but we prefer to give a self-contained
proof. We want to prove that $P^{-1} \in \mW^{0}$. In any case, there
exists a sequence of bounded operators $Q_n \in \mW^{0}$ such that
$Q_n \to P^{-1}$ in norm. (This is seen as follows, the functional
calculus applied to the bounded, invertible operator $PP^{*}$ shows
that there exists a sequence of polynomials $f_n$ such that
$f_{n}(PP^{*}) \to (PP^{*})^{-1}$. Then we can take $Q_{n} :=
P^{*}f_{n}(PP^{*})$).

Let $Q$ be a parametrix for $P$. That is, let $Q \in \mW^{0}$ be such
that $R := PQ - I \in \mW^{-\infty}$. (We use here the ellipticity of
$P$.) Then
$$
	P^{-1} = Q - P^{-1} R = \lim (Q - Q_n R),
$$
in the norm of bounded operators (the only one used in this proof).
Then $Q' := Q - Q_n R$ is also a parametrix of $P$ and is moreover
invertible if we choose $n$ large enough. If we write $PQ' = I + R'$,
then $I + R'$ is invertible as a bounded operator on $\mH$, and hence
$(I + R')^{-1} - I \in \mW^{-\infty}$, by axiom $(\psi)$. Thus $P^{-1}
= Q'(I + R')^{-1} \in \mW^{0}$. The proof is now complete.
\end{proof}

\section{Special holomorphic families\label{Sec.SHF}}

Let $\mW$ be an extended Weyl algebra, and 
let $\Omega$ be an open subset of $\CC$.

We shall need the following result.

\begin{proposition}\label{prop.2many}\ 
We have the following continuity properties

(i)\ The map $z \mapsto r^{z} \in S^{m}_{1, 0} := S^{m}_{1, 0}(V)$ is
holomorphic for $\re(z) < m$.

(ii)\ The functions $z \mapsto q(r^z)T  \in
\mW^{-\infty}_{k+1}$, $z \mapsto Tq(r^z) \in
\mW^{-\infty}_{k+1}$, and $z \mapsto q(r^z)\xi \in \mH$ are
holomorphic in $z \in \CC$ for any fixed $T \in \mW^{-\infty}_{k}$ and
any fixed $\xi \in \mW^{-\infty}\mH$.
\end{proposition}

\begin{proof}\ 
Recall that $r \in S^1_{cl} $ is a fixed symbol satisfying $r \ge 1$
and $r = \| v \|$ if $\| v \| \ge 2$. Fix $m \in \RR$. To prove that
the map
$$
	z \mapsto r^z \in S^m_{1, 0}
$$
is holomorphic on $\{ \re(z) < m \}$, we shall check that the complex
derivative
$$
	\pa_{z} r^{z} = \lim_{w \to 0} \frac{r^{z + w} - r^{z}}{w} =
	r^{z} \log r \in S^{m}_{1, 0}
$$
exists for $\re(z) < m$. By Proposition \ref{prop.cont.m}, it is
enough to prove this for $m = 0$.

Let $S^0_{1,0}([1,\infty))$ be the space of order zero symbols on
$[1,\infty)$, endowed with the usual semi-norms $p_n(a) := \sup_{x,\
k\leq n} |(x\pa_x)^{k} a|$.  We shall first check that the function
$$
	\{ \re(z) < 0 \} \ni z \mapsto x^{z} \in S^{0}([0,\infty))
$$ 
is holomorphic. Now $x^k \pa_x^{k} x^{z} = P(z) x^{z}$, with $P$ a
polynomial. Taylor's formula then gives that
\begin{equation}\label{eq.2many}
	|P(z + w)x^{z + w} - P(z)x^{z} - w\pa_{z}(P(z)x^{z})| \le
	C|w|^2 x^{-\epsilon} (\log x)^2 \le CC_\epsilon |w|^{2}
\end{equation}
for $\re(z) < -2 \epsilon$ and $|w| < \epsilon$, where $C$ depends
only on the polynomial $P$ and $C_\epsilon = \sup_{x \ge 1}
x^{-\epsilon} (\log x)^2 < \infty$. Since $x^{k}\pa_{x}$ is a linear
operator, this is enough to prove that
$$
	p_n\Big( \frac{x^{z + w} - x^{z}}{w} - x^{z} \log r \Big ) \to
	0, \quad \text{as} \; \; |w| \to 0,
$$
for any $n$. Consequently, $z \to x^z \in S^0_{1,0}([1,\infty))$ is
holomorphic for $\re(z) < 0$.

(i) then follows from the continuity of the map
$$
	S^{m}_{1,0}([1,\infty)) \ni a \to a \circ r \in S^{m}_{1, 0}
	:= S^{m}_{1, 0}(V),
$$
(take $\epsilon = 1$ in Lemma \ref{lemma.symb}).

The holomorphy on $\CC$ of the maps $z \mapsto q(r^{z}) T \in
\mW^{-\infty}_{k+1}$ and $T \mapsto a q(r^{z})\in \mW^{-\infty}_{k+1}$
then follows from Axioms \eqref{axiom.H} and \eqref{axiom.product}.
In turn, this then gives right away that the map $z \mapsto q(r^{z})
T\xi$, for $T \in \mW^{-\infty}$ and $\xi \in \mH$, is holomorphic,
and hence that $z \mapsto q(r^{z}) \xi$ is holomorphic on $\CC$ for
$\xi \in \mW^{-\infty}\mH$.
\end{proof}

We now introduce special holomorphic families.

\begin{definition}\label{def.holom.family}\
A family of operators $A(z) \in M_N(\mW^{mz + d})$, $m, d \in \RR$, is
called a {\em special holomorphic family of order $mz + d$ on
$\Omega$} if there exist $k \in \NN$, $b \in \mO(\Omega,
M_N(S^d_{cl}))$, and $R \in \mO(\Omega, M_N(\mW^{-\infty}_{k}))$ such
that
\begin{equation}\label{eq.spec.hol}
	A(z) = q(r^{mz}b(z)) + R(z).
\end{equation} 
A {\em special holomorphic family} is a special holomorphic family of
order $mz + d$ in a neighborhood of $0$, for some $m,d\in \RR$.
\end{definition}

Let us make now the simple observation that Lemma \ref{lemma.sh} shows
that $\sh{A(\overline{z})}$ is a special holomorphic family whenever
$A(z)$ is a special holomorphic family.

A special holomorphic family can be described in several ways.  Recall
that $\mO(\Omega, \mW^{w}) = \cup_{k \ge 1} \mO(\Omega,
\mW^{w}_{k})$, for any $w \in \CC \cup \{\infty\}$.

\begin{proposition}\label{prop.c.s.f}\ (i) 
A function $A(z) \in M_N(\mW^{mz + d})$ is a special holomorphic
family of order $mz + d$ on a neighborhood of $0$ if, and only if,
there exist holomorphic functions $B_1(z) \in M_N(\mW^d)$ and $R_1(z)
\in M_N(\mW^{-\infty})$, defined in a possibly smaller
neighborhood of $0$, such that $A(z) = q(r^{mz}I_N)B_1(z) + R_1(z)$.

(ii) Let $T \in \mO(\Omega, M_N(\mW^{-\infty}))$ and $A(z)$ be a
special holomorphic family on $\Omega$. Then $z \mapsto A(z)T(z)$ and
$z \mapsto T(z)A(z)$ are in $\mO(\Omega, M_N(\mW^{-\infty}))$.
\end{proposition}

\begin{proof}\ 
Let $N = 1$ and $d = 0$, for simplicity.  Assume $A(z) = q(r^{mz}a(z))
+ R(z)$ is a holomorphic family, as in the Definition
\ref{def.holom.family}. Let $b(z)$ be as in Corollary
\ref{cor.elliptic}, namely,
$$
	q(r^{mz})q(r^{-mz}b(z)) - I \in \mW^{-\infty}.
$$
We define $B_1(z) = q(r^{-mz}b(z))q(r^{mz}a(z)) \in \mW^{0}_2$. Then
$B_1(z)$ and $R_1(z) : = q(r^{mz})B(z) - A(z) \in \mW^{-\infty}$ are
the desired holomorphic families.

Conversely, assume that $A(z) = q(r^{mz}) B_1(z) + R_1(z)$, with $B_1$
and $R_1$ as in the statement of the proposition. By Proposition
\ref{prop.h.a.c}, there exists $c \in \mO(\Omega, S_{cl}^0)$ such that
$B_1(z) - q(c(z)) \in \mW^{-\infty}$. Let
$$
	a(z) = r^{-mz}\mP(r^{mz}, c(z)).
$$
Then $a \in \mO(\Omega, S_{cl}^0 )$ and $q(r^{mz})q(c(z)) -
q(r^{mz}a(z)) \in \mW^{-\infty}_{1}$ by Axiom \eqref{axiom.P} and
Proposition \ref{prop.2many}. Moreover,
\begin{multline*}
	R(z) := q(r^{mz})B_1(z) + R_1(z) - q(r^{mz}a(z)) \\ = \big
	(q(r^{mz})B_1(z) - q(r^{mz})q(c(z)) \big) + \big
	(q(r^{mz})q(c(z)) - q(r^{mz}a(z)) \big) + R_1(z) \in
	\mW^{-\infty}_{k}.
\end{multline*}

We continue to assume $N = 1$ and $d = 0$.  (ii) follows from the
descriptions of a special holomorphic family in (i) and Proposition
\ref{prop.2many}(ii) as follows. Let us write $A(z) = q(r^{mz} a(z)) +
R(z)$ for some $a \in \mO(\Omega, S_{cl}^{0})$ and $R(z) \in
\mO(\Omega, \mW^{-\infty}_{l})$. Choose $k$ such that $T \in
\mO(\Omega, \mW^{-\infty}_{k})$.

Then by \ref{prop.2many}(i), $r^{mz}a(z) \in S^{N}_{1, 0}$ is a
holomorphic family for $\re(z) < N/m$. The continuity of the
multiplication map $S^{N}_{1, 0} \times \mW^{-\infty}_{k} \ni (a, T)
\mapsto q(a) T \in \mW^{-\infty}_{k + 1}$ then proves that $q(r^{mz}a(z))
T(z)$ is a holomorphic function in $\mO(\Omega, \mW^{-\infty}_{k +
1})$.  Since $\mW^{-\infty}$ is an $LF$--algebra, $R(z)T(z)$ is also
in $\mO(\Omega, \mW^{-\infty}_{k + l})$. This proves (ii) for
$A(z)T(z)$. The proof for $T(z)A(z)$ is the same.
\end{proof}

Let us observe that by taking adjoints (and by replacing $z$ with
$\overline{z}$, any family of the form $A(z) = B_1(z)q(r^{mz}I_N) +
R_1(z)$, with $B_1(z) \in M_N(\mW^d)$ and $R_1(z) \in
M_N(\mW^{-\infty})$, is also a special holomorphic family of order $mz
+ d$.

The following consequence of the above Proposition will be useful later
on.

\begin{corollary}\label{cor.n.cont}\ 
Let $A(z) \in \mW^{mz}$ define a special holomorphic family. Then
$\|A(z)\|$ is continuous on the set $\Omega \cap \{\re(mz) < 0\}$ and
is bounded on the compact subsets of $\Omega \cap \{\re(mz) \le 0\}$.
\end{corollary}

\begin{proof}\ Let us write $A(z) = q(a(z)r^{mz}) + R(z)$ with
$a \in \mO(\Omega, S^0_{cl})$ and $R \in \mO(\Omega, \mW^{-\infty})$,
by Proposition \ref{prop.c.s.f}. Then $a(z)r^{mz}$ defines a
continuous function $\{\re(mz) < 0\} \to S^0_{1, 0}$, and hence $\|q(
a(z)r^{mz} ) \|$ is continuous on $\Omega \cap \{\re(mz) < 0\}$ by
Proposition \ref{prop.cont}(i).  For the second part of the statement,
let us observe that the map $z \mapsto r^{z} \in S^{0}_{1, 0}$ is bounded on
$\{\re(mz) \le 0\}$ and the map $z \mapsto a(z) \in S^0_{cl}$ is continuous
on $\Omega$. The result then follows from Proposition
\ref{prop.cont}(i).
\end{proof}

The asymptotic completeness of Proposition \ref{prop.h.a.c} extends to
special holomorphic families right away, giving the following
technical result, which will however be crucial in constructing
complex powers.

\begin{corollary}\label{cor.h.a.c}\ 
Let $A_i \in M_N(\mW^{mz + d - i})$ be special holomorphic families of 
order $mz + d - i$, $i \in \ZZ_+$, defined on some open subset $\Omega \in \CC$.  
Then there exists a special holomorphic family $A(z)$ of order $mz + d$ on
$\Omega$ such that $A(z) - \sum_{i = 0}^{k - 1} A_i(z)\in M_N(\mW^{mz
+ d - k})$.
\end{corollary}

\begin{proof}\ The proof is the same for all values of $N$, so we shall assume
that $N = 1$, for simplicity. Let us write $A_i(z) = q(r^{mz}a_i(z)) +
R_i(z)$, with $a_i \in \mO(\Omega, S^{d-i}_{cl} )$. We can choose by
Proposition \ref{prop.h.a.c} $a \in \mO(\Omega, S^{d}_{cl} )$ to be an
asymptotic sum of the families $a_i$, $i = 0, 1, 2, \ldots$. Then we
can take $A(z) = q(r^{mz}a(z))$.
\end{proof}

We conclude this section with a result that will allow us to
construct resolvents.

\begin{proposition}\label{prop.inv}\ 
Assume that the extended Weyl algebra satisfies the Conditions
$(\sigma)$ and $(\psi)$, and let $A(z)$ be a special holomorphic
family defined in an open neighbhorhood of $0$ in $\CC$ such that
$A(0) = I$.  Then there exists a special holomorphic family $A_1(z)$,
defined in a possibly smaller neighborhood of $0$ such that
$$
	A(z) A_1(-z) = A_1(-z)A(z) = I.
$$
\end{proposition}

\begin{proof}\ 
Let $A(z) = q(r^{mz})B(z) + R(z) \in M_N(\mW^{mz + d})$. We shall
assume $N = m = 1$ and $d = 0$, for simplicity. We can assume $b(0)=1$
and $R(0)= 0$. Then $C(z) := q(r^{-z})A(z)$ is a holomorphic function
with values in $\mW^0$ (defined in a small neighborhood of
$0$). Moreover, $C(z)$ is elliptic. By the holomorphic asymptotic
completeness (Proposition \ref{prop.h.a.c}), we can find a holomorphic
function $C_1(z) \in \mW^0$ such that 
$$
	R(z) := C_1(z)C(z) - I \in \mW^{-\infty}.
$$ 
By Proposition \ref{prop.cont.m}, $R$ is also holomorphic in a
possibly smaller neighborhood of $0$. We can also assume that $C_1(0)
= I$, and hence $R(0) = 0$.  By Axiom $(\psi)$, $(I + R(z))^{-1} = (I
+ R_1(z))$, for some holomorphic function $R_1(z)$ that is again
defined in a neighborhood of $0$. (This uses the fact that in a
$\Psi^*$-algebra the set of invertible elements is open and hence
inversion is continuous, by Banach's theorem mentioned also above, see
\cite{ban48}.) let
$$
	A_1(z) = (I + R_1(-z))C_1(-z)q(r^{z}).
$$ 
Then $A_1(-z)A(z) = I$. By Proposition \ref{prop.c.s.f}(i), $A_1(z)$ is a
special holomorphic family.
\end{proof}

\section{Complex powers of elliptic operators\label{Sec.CP}}

We now turn to the construction of complex powers of an elliptic,
strictly positive operator $A \in M_N(\mW^{m})$, $m \ge 0$, so that
$A^{m} \in M_N(\mW^{mz})$. We start by using Guillemin's method of
constructing complex powers in our slightly modified setting of
extended Weyl algebras. Then we adapt some arguments from
\cite{ShubinBook}.

We want to prove the following two theorems, analogues of Theorems 5.1
and 5.2 in \cite{Guil85}. For the ring of pseudodifferential
operators, these theorems are due to Seeley \cite{Seeley}. See also
\cite{ShubinBook}.

Let $\mW$ be an extended Weyl algebra and let $r \in S_{cl}^1 $ be as
before, namely $r \ge 1$, $r(\xi) = \|\xi\|$, for all $\xi \in V$,
such that $\|\xi\| \ge 2$. We shall assume that our extended Weyl
algebra satisfies condition $(\sigma)$.

\subsection{Symbolic complex powers}\
We first establish the following preliminary result, which is, in a
certain sense, a result about the symbolic part of complex powers.

\begin{proposition}\label{prop.interm}\ 
Recall that $\mW$ is an extended Weyl algebra satisfying Condition
$(\sigma)$.  Let $m \ge 0$, $A \in M_N(\mW^m)$ be elliptic, and
suppose that there exists $a \in \mO(\CC,M_N(S_{cl}^0 ))$, such that
\begin{equation}\label{eq.interm1}
	a(z+w) = a(z)a(w) \quad \text{and} \quad \sigma^{(m)}(A) =
	a(1)r^m.
\end{equation}
Then there exists a special holomorphic family $A(z) \in
M_N(\mW^{mz})$ such that
\begin{multline}\label{eq.interm2}
	\begin{gathered}
	A(z)A(w) - A(z + w) \in \mW^{-\infty},\quad
	\sigma^{(mz)}(A(z)) = a(z)r^{mz},\quad \text{and} \\ A(1) - A
	\in M_N(\mW^{-\infty}), \quad z, w \in \CC. 
	\end{gathered} 
\end{multline}
Moreover, this family is unique modulo smoothing operators.
\end{proposition}

Note that, by replacing $A(z)$ with $A(z) + I_{N} - A(0)$, we can
assume that $A(0) = I_{N} := \id_{\cH^{N}}$ in the Proposition above.

\begin{proof}\
The proof is done as in Guillemin's paper. We assume that the
reader is familiar with the details of that paper. Also, in the proof,
we shall assume $N = 1$, for simplicity of notation. (The extension to
the case $N > 1$ is non-trivial, however. To treat that case, we
proceed as in \cite{buci}.)

First we recall that the second holomorphic cohomology group of $\CC$
with values in $S^{\mu}_{cl}$ satisfies $H^{2}_{hol}(\CC,
S^{\mu}_{cl}) = 0$, \cite[Theorem 5.3]{Guil85}. For the benefit of the
reader, we now briefly recall the main steps of the proof.

The main idea is to approximate $A(z)$ by a sequence $A_n(z)$ such
that
$$
	A_n(z)A_n(z') - A_n(z + z') \in \mW^{m(z + z') - n - 1} \quad
	\text{and} \quad A_n(1) - A \in \mW^{m - n - 1}.
$$  
To start with, one can take $A_0(z) = q(a(z)r^{mz})$ (this works even
for $m = 0$).  Then one can solve $A_{n+1}(z) = A_n(z) + q(r^{mz}
b(z))$, where $b(z) \in S_{cl}^{- n -1}$, using the vanishing of
$H^2_{hol}(\CC, S^{\mu}_{cl})$. Moreover, at each step the solution is
unique, because two solutions differ by a linear map (modulo lower
order terms). This linear map is uniquely determined at by the
condition $A_n(1) - A \in \mW^{m -n -1}$.  The holomorphic asymptotic
completeness (Corollary \ref{cor.h.a.c}) then gives the result.
\end{proof}

\begin{theorem}\label{thmG5.1}\ 
Recall that $\mW$ is an extended Weyl algebra satisfying Conditions
$(\sigma)$ and $(\psi)$. Let $m \ge 0$, $A \in M_N(\mW^m)$, and $a(z)
\in M_N(S_{cl}^0) = M_N(S_{cl}^0(V))$, $z \in \CC$, satisfy
\eqref{eq.interm2} of the above Proposition.  Then we can find a
special holomorphic family $A_{z} \in M_N(\mW^{mz})$ such that
\begin{equation} 
	A_{z}A_{w}=A_{z + w}, \quad \sigma^{(mz)}(A_{z}) =
	a(z)r^{mz},\quad \text{and} \quad A_{1} - A \in
	M_N(\mW^{-\infty}), \quad .
\end{equation}
$z, w \in \CC$.
\end{theorem}

To prove Theorem \ref{thmG5.1}, we again proceed as in \cite{Guil85}.
This is contained in the following lemma.

Let us fix some notation before proceeding with the proof. We begin by
choosing $A(z)$ as in Proposition \ref{prop.interm}. We may assume
also that $A(z) = I_{N}$.  Then we write $A(z)A(w) = A(z + w) + F(z,
w)$, with $F: \CC^2 \to \mW^{-\infty}$ a holomorphic function. Let
$$
	P := [\pa_{z} A(z)]\vert_{z=0} \quad \text{and} \quad \text
	Q(z) := [\pa_{w} F(z, w)]\vert_{w = 0}.
$$
Then 
\begin{equation}\label{eq.diff}
	P \in q(S_{1, 0}^\epsilon) + \mW^{-\infty}, \;\; \epsilon > 0,
	\quad \text{and} \quad \pa_{z} A(z) = A(z) P - Q(z).
\end{equation}

\begin{lemma}\label{lemma.diff.eq}\ 
Let $P := [\pa_{z} A(z)]\vert_{z=0} \in q(S_{1, 0}^\epsilon) +
\mW^{-\infty},$ $\epsilon > 0,$ be as above.  Then there exists a
unique special holomorphic family of operators $A(z) \in \mW^{mz}$,
defined on $\CC$, such that $A_0 = I$,
$$
	\pa_{z} A_{z} = A_{z} P\quad \text{and} \quad A_{z} - A(z) \in
	\mW^{-\infty}.
$$ 
This family then also satisfies $A_{z}A_{w} = A_{z+w}$ and $A_{z}P =
PA_{z}$, for any $z, w \in \CC$.
\end{lemma}

\begin{proof}\ 
Let $A_1(z)$ be a special holomorphic family, defined in a small
neighborhood of $0$ in $\CC$ such that $A_1(-z) = A(z)^{-1}$, for all
$z$ close to $0$. The existence of this family is given by
Proposition~\ref{prop.inv} above. Then the function $z \mapsto
Q(z)A(z^{-1}) = Q(z)A_1(-z) \in \mW^{-\infty}$ is holomorphic, by
Proposition \ref{prop.c.s.f}(ii).

Let $\|\;\cdot\;\|_n$ be the sequence of Banach algebra norms on
$\mW^{-\infty}$ defining the topology on $\mW^{-\infty}$, as in
Condition $(\psi)$.  Denote by $X_n$ be the completion of
$\mW^{-\infty}$ in the norm $\|\;\cdot\;\|_n$, and let $R(z) \in X_n$,
with $z$ in a small neighborhood of $0$ in $\CC$ satisfy the {\em
holomorphic} differential equation
\begin{equation}
	\pa_{z}R(z) = (I + R(z)) Q(z)A(z)^{-1}.
\end{equation}
with initial condition $R(0) = 0$. Because this equation is
linear, a solution will exist on the domain $\Omega$ of
$Q(z)A(z)^{-1}$, by the existence theorem for ordinary differential
equations in Banach spaces. This solution does not depend on choice of
the norm $\|\;\cdot\;\|_{n}$, and hence $R\in \mW^{-\infty}$. 
%
Then $A_z := (I + R(z)) A(z)$ satisfies the differential equation $A_z' =
A_zP$, for $z \in \Omega$, as operators on $\mW^{-\infty}\mH$.  This
solution is unique.

Similarly, the solution of the differential equation $A(z)' = A(z)P$
with initial condition $A(0) = B$ is unique and is given by $A(z) = B
A_z$. This can be used to prove that
\begin{equation}\label{eq.mult}
	A_{z}A_{z'} = A_{z+z'},
\end{equation}
by differentiating with respect to $z'$. (Here $|z|, |z'| < \epsilon$,
for some $\epsilon > 0$ small enough such that the ball of radius
$2\epsilon$ is contained in $\Omega$.)  We define then $A_z =
(A_{z/n})^{n}$, where $n > |z|/\epsilon$.

By differentiating Equation \eqref{eq.mult} with respect to $z$, we
obtain the relation $A_zP = PA_z$, for any $z \in \CC$.
\end{proof}

\subsection{Semiclassical estimates} 
We have now proved Theorem~\ref{thmG5.1}. Before proceeding to
construct the complex powers, we need to establish some asymptotic
formulas, in the spirit of the estimates in \cite{ShubinBook} for
parameter dependent pseudodifferential operators. Recall that we
assume in this section that $\mW$ is an extended Weyl algebra
satisfying Condition $(\sigma)$.

Let $\mP_0(a,b) := q(a)q(b) - q(ab)$. Below we use crucially that
\begin{equation}\label{eq:mP_0-bd}
	\mP_0 : S^m_{1,0} \times S^{m'}_{1,0} \to \mW^{m+m'-1}
\end{equation}
is continuous (Axiom \ref{axiom.P}).  We fix an elliptic real-valued symbol $c\in S^m_{1,
0}$, $m>0$. (We note however that our arguments remain valid if 
$\arg
c $ is bounded away from $\alpha$, for some $\alpha \in [0,
2\pi]$.)  We shall need the following lemma, which is a particular
case of \cite[Lemma~2]{HV}. We include a proof, for completeness.

\begin{lemma}\label{lemma.bdd.a}\
The family $t\mapsto (c+it)^{-1}$, $|t|\geq 1$, is bounded in
$S^{-m}_{1, 0}$ and, moreover, $(c+it)^{-1} \to 0$ in $S^{-m +
\epsilon}_{1, 0}$, for any $\epsilon > 0$.
\end{lemma}

\begin{proof}\    
First, using Lemma \ref{lemma.symb}, we reduce to the case $m = 1$ and
$c = x \in S^{1}_{1,0}(\RR)$.

Then we notice that 
\begin{equation}\label{eq.ineq.a}
	|(x + \imath t)|^{-1} \leq C\<x\>^{-a} |t|^{a - 1}
\end{equation}
for any $t \ge 1$,  $0 \le a \le 1$, and $x \in \RR$, where $C$ is a
constant independent of $t$, $x$, or $a$.
Since $\pa_{x}^{k} (x + \imath t)^{-1} = (-1)^k k! (x + \imath t)^{-k
- 1}$, we obtain that $|x^{1 + k}\pa_{x}^{k} (x + \imath t)^{-1}|$ is
bounded, by using \eqref{eq.ineq.a} for $a = 1$.

To check the last part, we can assume $\epsilon < 1$. Using
\eqref{eq.ineq.a} for $a = 1 - \epsilon$, we obtain
$$
	|x^{1 - \epsilon + k}\pa_{x}^{k} (x + \imath t)^{-1}| \le Ck!
	|t|^{-\epsilon}
$$
for $|t| \ge 1$. This completes the proof.
\end{proof}

The following is an immediate consequence of the above lemma.

\begin{corollary}\label{cor.bdd2}\
The family $q(c+it)q((c+it)^{-1})-I$ is bounded in $\mW^{-1}$ for
$|t|\geq 1$. The proof that it converges to $0$ in $\mW^{-1 +
\epsilon}$, for any $\epsilon > 0$, is completely similar.
\end{corollary}

\begin{proof}\
First, $q((c+it)^{-1})$ is bounded in $\mW^{-m}$, by lemma
\ref{lemma.bdd.a}. Since $q(c)\in \mW^m$, we obtain from Proposition
\ref{prop.cont}(ii) that
\begin{equation*}
	q(c)q((c+it)^{-1}) - q(c(c+it)^{-1})
\end{equation*}
is bounded in $\mW^{-1}$ by \eqref{eq:mP_0-bd}. Thus,
\begin{equation*}\begin{split}
	q & (c+it)q((c+it)^{-1}) - I\\ & = \left[ q(c)q((c+it)^{-1}) -
	q(c(c+it)^{-1}) \right] + \left[
	q(c(c+it)^{-1})+itq((c+it)^{-1})\right] - I\\ & =\left[
	q(c)q((c+it)^{-1}) - q(c(c+it)^{-1})\right] +
	q((c+it)(c+it)^{-1}) - q(1) \\ & = q(c)q((c+it)^{-1}) -
	q(c(c+it)^{-1})
\end{split}\end{equation*}
is bounded in $\mW^{-1}$ and converges to $0$ in $\mW^{-1+\epsilon}$,
for any $\epsilon > 0$.
\end{proof}

We are now ready to prove one of our main results on resolvents in the
framework of extended Weyl algebras. Recall that in this section we
are assuming $\mW$ to satisfy condition $(\sigma)$.

\begin{theorem}\label{theorem.resolvent}\
Assume $\mW$ satisfies Condition $(\psi)$.  Let $T = q(c) + R$, with
$R \in \mW^{-\infty}$ and $c \in S^{m}_{1,0}$ real valued, elliptic,
as above. Then 
$$
	(T + \imath t)^{-1} = q((c + \imath t)^{-1}) (I + R_1(t)),
$$ 
where $R_1(t)$ is bounded in $\mW^{-1}$ and converges to $0$ in
$\mW^{0}$, as $|t| \to \infty$.

In particular, $(T + \imath t)^{-1}$ is bounded in $\mW^{-m}$ for $|t|
\ge 1$.
\end{theorem}

\begin{proof}\
Let us write $T = q(c) + R$, with $R \in \mW^{-\infty}$ and $c \in
S^{m}_{1,0}$ real valued, elliptic, as above. Then
$$
	R_2(t) := q((c + \imath t)^{-1}) (T + \imath t) - I = q((c +
	\imath t)^{-1})q(c + \imath t) + q((c + \imath t)^{-1})R - I
	\to 0
$$
in $\mW^{-1 + \epsilon}$, for any $\epsilon > 0$. In particular, if
$\epsilon = 1$, we obtain that $\|R_2(t)\| \to 0$ as $|t| \to \infty$,
by Proposition \ref{prop.cont}(ii). Hence $I + R_2(t)$ will be
invertible on $\mH$ for $|t|$ large. Assumption $(\psi)$ together with
Theorem \ref{theorem.invert} show that the set of invertible elements
of the Fr\'echet algebra $\mW^{0}$ is open. Hence, by Banach's theorem
mentioned also earlier, inversion is continuous on $\mW^{0}$. Using
the spectral invariance of $\mW^{0}$ (again by Theorem
\ref{theorem.invert}) and the continuity of the inversion for the
Fr\'echet algebra $\mW^{0}$, we obtain that $(I + R_2(t))^{-1} = I +
R_1(t)$, with $R_1(t) \in \mW^{0}$ and $R_1(t) \to 0$ in
$\mW^{0}$. This completes the proof.
\end{proof}

\begin{corollary}\label{cor.interest}\ 
Let $T \in \mW^{m}$, $m > 0$, $\sh{T} = T$, elliptic.  Then there
exists a parametrix family $t\mapsto G_t$, such that $G_t$ is bounded
in $\mW^{-m}$, converges to $0$ in $\mW^{-m + \epsilon}$, for any
$\epsilon > 0$, and $(T + \imath t)G_t - I$ is bounded in
$\mW^{-\infty}$.
\end{corollary}

\begin{proof}
Let $E_1(t)=I-q(c+it)q((c+it)^{-1}) - R$, so $E_1(t)$ is bounded in
$\mW^{-1}$.  Then the Neumann series $\sum_j E_1(t)^j$ can be summed
asymptotically, uniformly in $t$, to get a bounded family $E(t) \in
\mW^{-1}$. Thus $F_N(t) = E(t) - \sum_{j\leq N-1} E_1(t)^j$ is
uniformly bounded in $\mW^{-N}$, so
\begin{equation*}\begin{split}
	&(I-E_1(t))(I+E(t)) - I \\ & \qquad = (I-E_1(t))(I+\sum_{j\leq N-1}
	E_1(t)^j+F_N(t)) = E_1(t)^N+(I-E_1(t))F_N(t)
\end{split}\end{equation*}
is bounded in $\mW^{-N}$ as well.  Let $G_t=q((c+it)^{-1})(I+E(t))$.
Then
\begin{equation*}
	(q(c)+it)G_t-I=(I-E_1(t))(I+E(t))-I
\end{equation*}
is bounded in $\mW^{-N}$ for all $N$, hence in $\mW^{-\infty}$.

The convergence to $0$ is proved similarly.
\end{proof}

Of course, a similar calculation gives a left parametrix, and shows that
these may be taken equal. The parametrix identity 
then yields the following.

\begin{proposition}\label{prop.interest}
There exists a bounded family of operators $P_t\in\mW^{-m}$ such that
$(q(c)+it)^{-1}-P_t$ is uniformly bounded as a map between any
Sobolev spaces, in fact
\begin{equation*}
	\|(q(c)+it)^{-1}-P_t\|_{\mL(H^{(l)},H^{(l')})}\leq C|t|^{-1}.
\end{equation*}
\end{proposition}

\begin{proof}
We know that $(q(c)+it)G_t-I=E_t$ and $G_t(q(c)+it)-I=F_t$ are bounded
in $\mW^{-\infty}$. Thus,
\begin{equation*}\begin{split}
	(q(c)+it)^{-1}&=(G_t(q(c)+it)-F_t)(q(c)+it)^{-1}=G_t-F_t(q(c)+it)^{-1}\\
	&=G_t-F_t(q(c)+it)^{-1}((q(c)+it)G_t-E_t)\\
	&=G_t-F_tG_t+F_t(q(c)+it)^{-1}E_t.
\end{split}\end{equation*}
The first two terms are bounded in $\mW^{-m}$, respectively,
$\mW^{-\infty}$.  Moreover,
\begin{equation*}
	\|(q(c)+it)^{-1}\|_{\mL(L^2,L^2)}\leq |t|^{-1},
\end{equation*}
in particular bounded, and $E_t\in\mL(H^{(l)},\mH)$,
$F_t\in\mL(\mH,H^{(l')})$ are bounded, so the last term has the
desired properties.
\end{proof}

\subsection{Complex powers}
We shall now use this theorem to construct complex powers of positive
definite operators in $\mW$.

Recall that a (possibly unbounded) operator on a Hilbert space $\mH$
is {\em positive} if $(A\xi,\xi) \ge 0$ for any $\xi$ in the domain of
$A$, and is called {\em positive definite} if $(A\xi,\xi) \ge \epsilon
(\xi,\xi)$, for some $\epsilon > 0$, independent of $\xi$ in the
domain of $A$. If $A \in M_N(\mW^m)$ is elliptic and positive, then we
can define the complex powers $A^z$ by the spectral theorem, as
above. To prove that $A^z$ are pseudodifferential operators, we assume
below that $A$ is also invertible (and hence positive definite).
(Recall that every positive $A \in \mW$ is automatically
essentially self-adjoint by Proposition \ref{prop.cont}.)

We are now ready to state and prove the main result of this section.

\begin{theorem}\label{thmG5.2}\ 
Let $\mW$ be an extended Weyl algebra satisfying Conditions $(\sigma)$
and $(\psi)$. Let $m \ge 0$ and $A \in M_N(\mW^m)$ be positive
definite and elliptic. 
Assume $\sigma^{(m)}(A) > 0$. Then $A^{z}\in
M_N(\mW^{mz})$, for all $z \in \CC$, and defines a special holomorphic
family such $\sigma^{(mz)}(A^{z}) = \sigma^{(m)}(A)^{z}$.

In particular, $A^{z} = q(a(z)r^{mz}) + R(z)$, with $a(z) \in \mO(\CC,
S^0_{cl})$ and $R(z) \in \mO(\CC, \mW^{-\infty})$.
\end{theorem}

\begin{proof}\ 
The last part is nothing but the definition of a special holomorphic
family (Definition \ref{def.holom.family} and Proposition
\ref{prop.c.s.f}). We shall assume $N = 1$, for simplicity. Recall
that $\sh{T} := T^*\vert_{\mW^{-\infty}\mH}$ (Lemma
\ref{lemma.sh}). 

The positivity of $A$ shows that $\sigma^{(m)}(\xi) > 0$ if $\xi \not
= 0$.  We can then define the powers $\sigma^{(m)}(A)^{z}$. Then we
apply Proposition \ref{prop.interm} to $a(z) = \sigma^{(m)}(A)^{z}$,
obtaining the family $A(z)$.  By replacing $A(z)$ with $1/2(A(z) +
\sh{A(\overline{z})})$, if necessary, we can also assume
$A(\overline{z}) = \sh{A(z)}$.

Then we construct $A_z$ as in the proof of Lemma
\ref{lemma.diff.eq}. The family $\sh{A_{\overline z}}$ will satisfy
the same differential equation as $A_z$ (namely $\pa_{z}A_z = A_zP = P
A_{z}$, using the notation of Lemma \ref{lemma.diff.eq}). By the
uniqueness of the solutions of this equation with $A_0 = I$, we obtain
\begin{equation}\label{eq.s.a}
	A_{\overline z} = \sh{A_z}.
\end{equation}
The resulting family $A_z \in \mW^{mz}$ is then a special holomorphic
family satisfying
\begin{equation*}
	A_0 = I, \quad A_1 - A \in \mW^{-\infty}, \quad
	\sigma^{(mz)}(A_z) = a(z), \quad \text{and} \quad A_{z + w} =
	A_{z} A_{w}, \quad z, w \in \CC.
\end{equation*}
It satisfies the additional condition $A_{\imath s}^* = A_{-\imath
s}$, for $s \in \RR$, by Equation \eqref{eq.s.a}. Since the function
$A_z\xi$ is holomorphic for $\xi \in \mW^{-\infty}\mH$, we obtain that
$A_{\imath s}$, $s \in \RR$, is a strongly continuous one-parameter
group of unitaries.

Let us prove now that $A_z = (A_{1})^{z}$. By Stone's theorem
\cite{ReedSimon}, we have that $A_{\imath s} = e^{\imath sT}$, where
$T$ is a self-adjoint operator and $\pa_s A_{\imath s} \xi = A_{\imath
s}T \xi$, for any $\xi$ in the domain of $T$ (this is the actual
definition of the domain of $T$). In particular, the domain of $T$
contains $\mW^{-\infty}\mH$ (because $A_{is}\xi$ is differentiable for
$\xi \in \mW^{-\infty}\mH$) and $T = P$ on $\mW^{-\infty}\mH$, by the
Cauchy-Riemann equations.

We claim that $\sigma(T) = \sigma(P) \subset [-M, \infty)$. Indeed, we
know from Corollary \ref{cor.n.cont} that the function $z \mapsto \|A_z\|$
is bounded on the compact subsets of $\re(z) \le 0$. Let $\|A_z\| \le
M$, if $|z| \le 1$ and $\re(z) \le 0$. Then $\|A_z\| \le M^{|z| + 1}$
for all $\re(z) \le 0$. Let
\begin{equation}
	R(\lambda):= \int_{0}^\infty e^{\lambda t}A_{-t}dt,
\end{equation}
which is convergent if $\re(\lambda) < - M$ and satisfies $(P -
\lambda) R(\lambda) = R(\lambda) (P - \lambda) = I$.

We now check that $A_{\imath s + t}\xi = e^{(\imath s + t)T}\xi$, for
any $s, t \in \RR$ and any $\xi \in \mW^{-\infty}\mH$. We have already
checked that this is true for $t = 0$. Fix $s$ and $\xi$. Because
$\sigma(T) \subset [-M, \infty)$, $e^{(\imath s + t)T}$ is bounded for
$t \le 0$, and hence the function
$$
	F(t) := A_{\imath s + t}\xi - e^{(\imath s + t)T} \xi
$$
is defined and continuous for $t \le 0$. By the same argument, $F$ is
differentiable when $t < 0$. Its derivative is $F'(t) = 0$, for $t <
0$, because $T\xi = P\xi$, $\xi \in \mW^{-\infty}\mH$. This proves
that $F(t) = 0$ for all $t \le 0$. Therefore $A_{\imath s + t}\xi =
e^{(\imath s + t)T}\xi$, for any $s, t \in \RR$, $t \le 0$, and any
$\xi \in \mW^{-\infty}\mH$.  For $t > 0$, we have
$$
	A_{\imath s + t}\xi = (A_{-\imath s - t})^{-1}\xi =
	(e^{-(\imath s + t)T})^{-1}\xi = e^{(\imath s + t)T}\xi,
$$
for any $\xi \in \mW^{-\infty}\mH$.

Finally, 
$$
	A_z = e^{\imath zT} = (e^{\imath T})^{z} = (A_1)^{z}.
$$

We now use the results of the previous subsection. The resolvent
identity gives for $\re(z) <-1$ that
$$
	A^{z}= (A_{1})^{z} + \left(A^{z}-(A_{1})^{z}\right) = A_{z} +
	\frac{1}{2\pi i} \int_{\gamma}\lambda^{z}(\lambda I -
	A)^{-1}(A_{1}-A) (\lambda I - A_{1})^{-1} d\lambda.
$$
We can now use then Theorem \ref{theorem.resolvent} to conclude that
the right hand side of the above integral converges in
$\mW^{-\infty}$.

The full conclusion follows using $A^{z+k}=A^k A^z$ with $\re z<-1$.
\end{proof}

This yields the following results (due to Seeley \cite{Seeley}, Schrohe
\cite{Schrohe}, and Kordyukov \cite{kordyukov2}).

\begin{corollary}\ 
Let $A \in \Psi^{m}(M)$ (respectively, $A \in \Psi^{m}_{sc}(\RR^n)$,
or $A \in B\Psi^m(M)$) (see Examples 1-3). Assume $A > 0$, $m \ge 0$,
and $A$ elliptic if $m > 0$. Then $A^{z} \in \Psi^{m}(M)$
(respectively, $A \in \Psi^{z}_{sc}(\RR^n)$, or $A \in B\Psi^{z}(M)$).
\end{corollary}

\begin{appendix}
\section{Semi-ideals and holomorphic functional calculus}
\label{secb}

For the convenience of the reader, we include a special case of a
construction in \cite{lmn2} that was used in Section \ref{Sec.AS} to
define an algebra that is closed under holomorphic functional calculus
containing a given Weyl algebra.

Recall that a subspace $\cJ\subseteq \cB$ of a unital algebra $\cB$ is
said to be a {\em semi-ideal} provided we have $xby\in \cJ$ for all
$x,y\in\cJ$ and all $b\in \cB$.  Suppose in addition that $\cB$ is a
unital $C^{*}$-algebra. Then we have $f(x)\in\cJ$ for all $f$ that are
holomorphic in a neighborhood of the spectrum $\sigma_{\cB}(x)$ of
$x\in\cJ$ and satisfy $f(0)=0$, i.e.\ the semi-ideal $\cJ$ is closed
under holomorphic functional calculus in $\cB$.

Let now $T:\cH\supseteq\cD(T)\rightarrow\cH$ be a closed operator and
$\fA \subset \mL(\mH)$ be a closed, symmetric subalgebra (\ie a
$C^*$--subalgebra). We are going to associate a family of semi-ideals
$\cJ_{k}:=\cJ_{k}(T)$, $k\ge 0$, to $T$.  Let $\cJ_{0}:=\fA\subset
\cL(\cH)$ be a $C^{*}$-subalgebra, and let $\cJ_{1}$ be the space of
all $A\in\cJ_{0}$ such that
\begin{enumerate}
\item $A(\cH)\subseteq\cD(T)$ and $\otl(A):=TA\in\cJ_{0}$.
\item There exists $\otr(A)\in\cJ_{0}$ with 
      $\otr(A)f=ATf$ for all $f\in\cD(T)$. 
\item $A(\cH)\subseteq\cD(T)$, and there exists 
      $\otlr(A)\in\cJ_{0}$ such that
      $\otlr(A)f=TATf$ for all $f\in\cD(T)$.
\end{enumerate} 
Moreover, let $\cJ_{k+1}$ be the space of all $A\in\cJ_{k}$ with
$\otl(A),\otr(A),\otlr(A)\in\cJ_{k}$.  The spaces $\cJ_{k}$ are
naturally endowed with the following norms: for $A\in\cJ_{0}$, let
$p_{0}(A)=\norm{A}{\cB(\cH)}$, and for $ A\in\cJ_{k+1}$ let
$$p_{k+1}(A):=p_{k}(A)+p_{k}(\otl(A))+p_{k}(\otr(A))
+p_{k}(\otlr(A))\,.$$ The projective limit
$\cJi:=\bigcap_{k\ge 0}\cJ_{k}$ is given the topology induced by the
system of norms $(p_{k})_{k\ge 0}$.  The main properties of this
construction are summarized in the following Proposition; the
straightforward proof can be found in \cite{lmn2, LN0}.

\begin{proposition}
\label{semi}

\noindent
\begin{enumerate}
\item For $k\ge 0$, $(\cJ_{k},p_{k})$ is a unital Banach algebra.
\item $(\cJi,(p_{k})_{k\ge 0})$ is a submultiplicative
      Fr\'echet algebra.
\item The canonical inclusion 
      $\cJ_{k}\hookrightarrow\cJ_{0}$ is continuous for all
      $k\ge 0\cup\{\infty\}$.
\item For $k\ge 0\cup\{\infty\}$, 
      $\cJ_{k}$ is a semi-ideal in the $C^{*}$-algebra $\fA$; in particular,
      $\cJ_{k}$ is closed under holomorphic functional calculus in
      $\cB(\cH)$. Moreover, the canonical map
      $$\cJ_{k}\times\fA\times\cJ_{k}\rightarrow\cJ_{k}$$
      is jointly continuous.
\end{enumerate}
\end{proposition}

\begin{remark}
\label{remstar}
Note that it is not clear, and in general not true that the
spaces $\cJ_{k}$ are symmetric subspaces of $\cB(\cH)$. However,
we easily obtain this property by considering the spaces
$$\cJs_{k}:=\{A\in\cJ_{k}:A^{*}\in\cJ_{k}\}\,, k\ge
0\cup\{\infty\}\,.$$ Note that Proposition \ref{semi} remains true for
these smaller spaces.
\end{remark}

\begin{remark}
\label{rough}
Roughly speaking, $A\in\fA$ is in $\cJ_{k}$ if for all $0\leq N,N'\leq
k$, the operators $T^{N}AT^{N'}$ extend to bounded operators on $\cH$
that belong to the given $C^{*}$-subalgebra $\fA\subseteq\cB(\cH)$.
\end{remark}

Furthermore, note that naturally associated to
$T:\cH\supseteq\cD(T)\rightarrow\cH$ there is a discrete scale of
Hilbert spaces, namely
\begin{eqnarray*}
\hTn & :=  & {\cH}\\
\hTa{1} & := & \cD(T)\\
\hTk & := & \{f\in \hTa{k-1}:Tf\in\hTa{k-1}\}\,, k\geq2\,.
\end{eqnarray*}
The spaces $\hTk$ are endowed with the iterated graph norms with
respect to $T$, i.e.\ $q(f):=\norm{f}{\cH}$ and
$q_{k}(f):=q_{k-1}(f)+q_{k}(Tf)$, $f\in\hTk$, $k\geq1$.  The Hilbert
space $\hTk$ is said to be the {\em Sobolev space of order $k$
associated to $T$}.  The intersection
$\hTu:=\bigcap_{k=0}^{\infty}\hTk$ is endowed with the system of norms
$(q_{k})_{k\ge 0}$.

\begin{proposition}
\label{propa4}
For $k\ge 0\cup\{\infty\}$, the canonical bilinear map
$$\cJ_{k}\times\hTn\longrightarrow \hTk:
(A,f)\longmapsto Af$$
is well-defined and continuous.
\end{proposition}

\begin{proof}
This follows immediately by induction from the definitions.
\end{proof}
\end{appendix}


\end{document}